\documentclass[11pt]{amsart}
\usepackage{booktabs}
\usepackage[section]{placeins}
\usepackage{bmpsize}
\usepackage{graphicx}
\usepackage{amssymb,amsmath}
\usepackage{caption,subcaption}
\usepackage{color}
\usepackage{bm}
\usepackage{wrapfig}
\usepackage{xcolor}
\usepackage[symbol]{footmisc}
\usepackage{afterpage}
\usepackage{cite}
\usepackage{hyperref}

\usepackage[linesnumbered,ruled,vlined]{algorithm2e}

\SetCommentSty{mycommfont}
\usepackage{amsthm}%
\usepackage{thmtools,thm-restate}
\usepackage[top=25mm, bottom=25mm, left=30mm, right=30mm]{geometry}

\usepackage{multirow}%
\usepackage{amsmath,amssymb,amsfonts}%
\usepackage{mathrsfs}%
\usepackage[title]{appendix}%
\usepackage{xcolor}%
\usepackage{textcomp}%
\usepackage{manyfoot}%
\usepackage{booktabs}%
\usepackage{listings}%
\usepackage{nicematrix}

\usepackage{tikz}
\usepackage{caption}
\usetikzlibrary{decorations.pathreplacing,calligraphy}
\usetikzlibrary{calc,patterns}

\usetikzlibrary{fit}

\newtheorem{remark}{Remark}%

\newcommand{\lsp}{\vspace{3mm}}

\newcommand{\vct}[1]{\bm{\mathsf{#1}}}
\newcommand{\mtx}[1]{\bm{\mathsf{#1}}}

\newcommand{\tikzcuboid}[5]{% width, height, depth, scale, buffer
\begin{tikzpicture}[scale=#4]
\foreach \x in {0,...,#1}
{   \draw (\x ,0  ,#3 ) -- (\x ,#2 ,#3 );
    \draw (\x ,#2 ,#3 ) -- (\x ,#2 ,0  );
}
\foreach \x in {0,...,#2}
{   \draw (#1 ,\x ,#3 ) -- (#1 ,\x ,0  );
    \draw (0  ,\x ,#3 ) -- (#1 ,\x ,#3 );
}
\foreach \x in {0,...,#3}
{   \draw (#1 ,0  ,\x ) -- (#1 ,#2 ,\x );
    \draw (0  ,#2 ,\x ) -- (#1 ,#2 ,\x );
}

\foreach \x in {0,...,#1}{
    \pgfmathtruncatemacro{\j}{mod(\x,#5)}
    \foreach \y in {0,...,#2}{
        \ifthenelse{\equal{\j}{0}}
        {\filldraw[red] (\x,\y,#3) circle (4pt);}
        {\filldraw[blue] (\x,\y,#3) circle (4pt);}
    }
}
\foreach \y in {0,...,#2}{
    \foreach \z in {0,...,#3}{
        \filldraw[red] (#1,\y,\z) circle (4pt);
    }
}
\foreach \x in {0,...,#1}{
    \pgfmathtruncatemacro{\j}{mod(\x,#5)}
    \foreach \z in {0,...,#3}{
        \ifthenelse{\equal{\j}{0}}
        {\filldraw[red] (\x,#2,\z) circle (4pt);}
        {\filldraw[blue] (\x,#2,\z) circle (4pt);}
    }
}
\pgfmathtruncatemacro{\numbuf}{#1/#5}
\pgfmathsetmacro{\halfbuf}{#5/2}
\foreach \x in {0,...,\numbuf}{
    \pgfmathtruncatemacro{\odd}{2*\x+1}
    \node [below,red] at (\x*#5,0,#3) {$I_{\odd}$};
}
\pgfmathtruncatemacro{\numbuf}{#1/#5 - 1}
\pgfmathsetmacro{\halfbuf}{#5/2}
\foreach \x in {0,...,\numbuf}{
    \pgfmathtruncatemacro{\even}{2*\x+2}
    \node [below,blue] at (\x*#5+\halfbuf,0,#3) {$I_{\even}$};
}

\draw [pen colour={black},
    decorate,
    thick,
    decoration = {calligraphic brace,
        raise=5pt,
        amplitude=5pt,
        aspect=0.5,mirror}] (#1,0) --  (#1,#2)
node[pos=0.5,right=10pt,black]{$n$};

\draw [pen colour={black},
    decorate,
    thick,
    decoration = {calligraphic brace,
        raise=5pt,
        amplitude=5pt,
        aspect=0.5}] (#1+0.5-#5,#2) --  (#1-0.5,#2)
node[pos=0.5,above=10pt,black]{$b$};
\end{tikzpicture}
}

\newcommand{\reducedtikzcuboid}[5]{% width, height, depth, scale, buffer
\begin{tikzpicture}[scale=#4]
\foreach \x in {0,...,#1}
{   \draw (\x ,0  ,#3 ) -- (\x ,#2 ,#3 );
    \draw (\x ,#2 ,#3 ) -- (\x ,#2 ,0  );
}
\foreach \x in {0,...,#2}
{   \draw (#1 ,\x ,#3 ) -- (#1 ,\x ,0  );
    \draw (0  ,\x ,#3 ) -- (#1 ,\x ,#3 );
}
\foreach \x in {0,...,#3}
{   \draw (#1 ,0  ,\x ) -- (#1 ,#2 ,\x );
    \draw (0  ,#2 ,\x ) -- (#1 ,#2 ,\x );
}

\foreach \x in {0,...,#1}{
    \pgfmathtruncatemacro{\j}{mod(\x,#5)}
    \foreach \y in {0,...,#2}{
        \ifthenelse{\equal{\j}{0}}
        {\filldraw[red] (\x,\y,#3) circle (4pt);}
        {\filldraw[gray] (\x,\y,#3) circle (4pt);}
    }
}
\foreach \y in {0,...,#2}{
    \foreach \z in {0,...,#3}{
        \filldraw[red] (#1,\y,\z) circle (4pt);
    }
}
\foreach \x in {0,...,#1}{
    \pgfmathtruncatemacro{\j}{mod(\x,#5)}
    \foreach \z in {0,...,#3}{
        \ifthenelse{\equal{\j}{0}}
        {\filldraw[red] (\x,#2,\z) circle (4pt);}
        {\filldraw[gray] (\x,#2,\z) circle (4pt);}
    }
}
\pgfmathtruncatemacro{\numbuf}{#1/#5}
\pgfmathsetmacro{\halfbuf}{#5/2}
\foreach \x in {0,...,\numbuf}{
    \pgfmathtruncatemacro{\odd}{2*\x+1}
    \node [below,red] at (\x*#5,0,#3) {$I_{\odd}$};
}

\draw [pen colour={black},
    decorate,
    thick,
    decoration = {calligraphic brace,
        raise=5pt,
        amplitude=5pt,
        aspect=0.5,mirror}] (#1,0) --  (#1,#2)
node[pos=0.5,right=10pt,black]{$n$};

\draw [pen colour={black},
    decorate,
    thick,
    decoration = {calligraphic brace,
        raise=5pt,
        amplitude=5pt,
        aspect=0.5}] (#1+0.5-#5,#2) --  (#1-0.5,#2)
node[pos=0.5,above=10pt,black]{$b$};
\end{tikzpicture}
}

\newcommand{\sparseA}[4]{% buffer, nsq, pans, scale
    \begin{tikzpicture}[scale=#4]

    \pgfmathsetmacro{\b}{#2}
    \pgfmathsetmacro{\nsq}{#1}
    \pgfmathsetmacro{\pans}{#3}
    \pgfmathsetmacro{\tot}{\nsq * (\pans + 1) + \b * \nsq * \pans}
        \draw[black](-\nsq,+\nsq) rectangle (\tot-\nsq,\tot+\nsq);
    \foreach \i in {1,...,\pans}{
        \pgfmathsetmacro{\startd}{(\i-1) * \b * \nsq + (\i-1) * \nsq}
        \pgfmathsetmacro{\endd}{\i * \b * \nsq + (\i-1) * \nsq}
        \pgfmathtruncatemacro{\l}{2*(\i)}
        \pgfmathtruncatemacro{\lpl}{2*(\i)+1}
        \pgfmathtruncatemacro{\lmn}{2*(\i)-1}
        \draw[blue] (\startd,\tot - \startd) rectangle (\endd,\tot - \endd) node[pos=.5] {$\mtx A_{\l\l}$}; % diagonal element
        \draw[blue] (\startd, \tot - \startd) rectangle (\endd, \tot - \startd + \nsq) node[pos=.5, scale=0.5] {$\mtx A_{\lmn\l}$}; % A_LC
        \draw[blue] (\startd, \tot - \startd - \nsq*\b - \nsq) rectangle (\endd, \tot - \startd + \nsq - \nsq*\b - \nsq) node[pos=.5,scale=0.5] {$\mtx A_{\lpl\l}$}; % A_RC
        \draw[blue] (\startd - \nsq, \tot - \startd) rectangle (\endd - \b * \nsq, \tot - \endd) node [pos=.5,scale=0.5] {$\mtx A_{\l\lmn}$}; % A_CL
        \draw[blue] (\startd  + \b * \nsq, \tot - \startd) rectangle (\endd + \nsq, \tot - \endd) node [pos=.5, scale=0.5] {$\mtx A_{\l\lpl}$}; % A_CL
    }
    \pgfmathsetmacro{\panspl}{\pans+1}
    \foreach \i in {1,...,\panspl}{
        \pgfmathsetmacro{\startd}{(\i-1) * \b * \nsq + (\i-2) * \nsq}
        \pgfmathsetmacro{\endd}{(\i-1) * \b * \nsq + (\i-1) * \nsq}
        \pgfmathtruncatemacro{\l}{2*(\i-1) + 1}
        \draw[red] (\startd,\tot- \startd) rectangle (\endd,\tot - \endd) node [pos = .5, scale=0.5] {$\mtx A_{\l\l}$}; % diagonal element
    }
    \pgfmathsetmacro{\tmp}{0.5*\b}
    \node [below,black] at (\tmp * \nsq * \pans+ \nsq,\nsq) {$N \times N$};
    \end{tikzpicture}
}

\newcommand{\denseT}[3]{%nsq,pans,scale
    \begin{tikzpicture}[scale=#3]
    \pgfmathsetmacro{\b}{1}
    \pgfmathsetmacro{\nsq}{#1}
    \pgfmathsetmacro{\pans}{#2}
    \pgfmathsetmacro{\tot}{\nsq * (\pans + 1) + \b * \nsq * \pans}
        \draw[black](-\nsq,+\nsq) rectangle (\tot-\nsq,\tot+\nsq);
    \foreach \i in {1,...,\pans}{
        \pgfmathsetmacro{\startd}{(\i-1) * \nsq + (\i-1) * \nsq}
        \pgfmathsetmacro{\endd}{\i * \nsq + (\i-1) * \nsq}
        \pgfmathtruncatemacro{\l}{4*(\i-1)+3}
        \pgfmathtruncatemacro{\lpl}{4*(\i)+1}
        \pgfmathtruncatemacro{\lmn}{4*(\i)-3}
        \draw[fill=red, fill opacity=0.5,text opacity=1] (\startd,\tot - \startd) rectangle (\endd,\tot - \endd) node[pos=.5,scale=0.6] {$\mtx T_{\l\l}$}; % diagonal element
        \draw[fill=red, fill opacity=0.5,text opacity=1] (\startd, \tot - \startd) rectangle (\endd, \tot - \startd + \nsq) node[pos=.5, scale=0.6] {$\mtx T_{\lmn\l}$}; % A_LC
        \draw[fill=red, fill opacity=0.5,text opacity=1] (\startd, \tot - \startd - \nsq*\b - \nsq) rectangle (\endd, \tot - \startd + \nsq - \nsq*\b - \nsq) node[pos=.5,scale=0.6] {$\mtx T_{\lpl\l}$}; % A_RC
        \draw[fill=red, fill opacity=0.5,text opacity=1] (\startd - \nsq, \tot - \startd) rectangle (\endd - \b * \nsq, \tot - \endd) node [pos=.5,scale=0.6] {$\mtx T_{\l\lmn}$}; % A_CL
        \draw[fill=red, fill opacity=0.5,text opacity=1] (\startd  + \b * \nsq, \tot - \startd) rectangle (\endd + \nsq, \tot - \endd) node [pos=.5, scale=0.6] {$\mtx T_{\l\lpl}$}; % A_CL
    }
    \pgfmathsetmacro{\panspl}{\pans+1}
    \foreach \i in {1,...,\panspl}{
        \pgfmathsetmacro{\startd}{(\i-1) * \b * \nsq + (\i-2) * \nsq}
        \pgfmathsetmacro{\endd}{(\i-1) * \b * \nsq + (\i-1) * \nsq}
        \pgfmathtruncatemacro{\l}{4*(\i-1) + 1}
        \draw[fill=red, fill opacity=0.5,text opacity=1] (\startd,\tot- \startd) rectangle (\endd,\tot - \endd) node [pos = .5, scale=0.6] {$\mtx T_{\l\l}$}; % diagonal element
    }
    \node [below,black] at (\nsq*\pans-\nsq/2,\nsq) {$\frac N b \times \frac N b$};
    \end{tikzpicture}
}

\newcommand{\simpleproofpicture}[3]{ %n,b,scale
\begin{tikzpicture}[scale=#3]
\pgfmathsetmacro{\b}{#2}
\pgfmathsetmacro{\n}{#1}

\foreach \x in {0,...,\b}
   \draw (\x ,0 ) -- (\x ,\n );

\foreach \x in {0,...,\n}
    \draw (0  ,\x ) -- (\b ,\x );

\foreach \x in {0,...,\b}
{
    \foreach \y in {0,...,\n}
    \filldraw[gray] (\x,\y) circle (4pt);
}

\foreach \y in {0,...,\n}
    \filldraw[black] (0,\y) circle (4pt);

\pgfmathsetmacro{\tmpstart}{2*\b}
\pgfmathsetmacro{\tmpend}{4*\b}
\foreach \y in {\tmpstart,...,\tmpend}
    \filldraw[red] (0,\y) circle (4pt);

\node [left,black] at (0,0) {$I_1$};
\draw [pen colour={black},
    decorate,
    thick,
    decoration = {calligraphic brace,
        raise=5pt,
        amplitude=5pt,
        aspect=0.5}] (0,2*\b) --  (0,4*\b)
node[pos=0.5,left=10pt,black]{$J_B$};

\draw [pen colour={black},
    decorate,
    thick,
    decoration = {calligraphic brace,
        raise=5pt,
        amplitude=5pt,
        aspect=0.5,mirror}] (0,0) --  (\b,0)
node[pos=0.5,below=10pt,black]{$b$};
\end{tikzpicture}
}

\newcommand{\detailedproofpicture}[3]{% n,b,scale
\begin{tikzpicture}[scale=#3]
\pgfmathsetmacro{\b}{#2}
\pgfmathsetmacro{\n}{#1}

\foreach \x in {0,...,\b}
   \draw (\x ,0 ) -- (\x ,\n );

\foreach \x in {0,...,\n}
    \draw (0  ,\x ) -- (\b ,\x );

\foreach \x in {1,...,\b}
{
    \foreach \y in {0,...,\n}
    \filldraw[gray] (\x,\y) circle (4pt);
}

\foreach \y in {0,...,\n}
    \filldraw[black] (0,\y) circle (4pt);

\pgfmathsetmacro{\tmpstart}{2*\b}
\pgfmathsetmacro{\tmpend}{4*\b}
\pgfmathsetmacro{\tmpn}{\n-1}
\foreach \y in {\tmpstart,...,\tmpend}
    \filldraw[red] (0,\y) circle (4pt);

\foreach \x in {1,...,\b}
{
    \filldraw[blue] (\x,\tmpstart) circle (4pt);
    \filldraw[blue] (\x,\tmpend) circle (4pt);
}

\node [left,black] at (0,0) {$I_1$};
\node [right,black] at (\b,\tmpstart) {$J_{\gamma}$};
\node [right,black] at (\b,\tmpend) {$J_{\gamma}$};

\draw [pen colour={black},
    decorate,
    thick,
    decoration = {calligraphic brace,
        raise=5pt,
        amplitude=5pt, mirror,
        aspect=0.5}] (\b,0) --  (\b,\tmpstart-1)
node[pos=0.5,right=10pt,black]{$J_{\alpha}$};

\draw [pen colour={black},
    decorate,
    thick,
    decoration = {calligraphic brace,
        raise=5pt,
        amplitude=5pt, mirror,
        aspect=0.5}] (\b,\tmpstart+1) --  (\b,\tmpend-1)
node[pos=0.5,right=10pt,black]{$J_{\beta}$};

\draw [pen colour={black},
    decorate,
    thick,
    decoration = {calligraphic brace,
        raise=5pt,
        amplitude=5pt, mirror,
        aspect=0.5}] (\b,\tmpend+1) --  (\b,\n)
node[pos=0.5,right=10pt,black]{$J_{\alpha}$};

\draw [pen colour={black},
    decorate,
    thick,
    decoration = {calligraphic brace,
        raise=5pt,
        amplitude=5pt,
        aspect=0.5}] (0,\tmpstart) --  (0,\tmpend)
node[pos=0.5,left=10pt,black]{$J_B$};

\draw [pen colour={black},
    decorate,
    thick,
    decoration = {calligraphic brace,
        raise=5pt,
        amplitude=5pt,
        aspect=0.5,mirror}] (0.5,0) --  (\b-0.5,0)
node[pos=0.5,below=10pt,black]{$b$};
\end{tikzpicture}
}

\begin{document}
\begin{center}
\textsc{SlabLU: A Two-Level Sparse Direct Solver for Elliptic PDEs}

\lsp

Anna Yesypenko\footnotemark[1] and Per-Gunnar Martinsson\footnotemark[2]

\footnotetext[1]{Oden Institute, University of Texas at Austin. Email: \texttt{annayesy@utexas.edu}}
\footnotetext[2]{Oden Institute, University of Texas at Austin. Email: \texttt{pgm@oden.utexas.edu}}

\lsp

\begin{minipage}{130mm}
\small

\textbf{Abstract:} The paper describes a sparse direct solver for the linear systems that arise from the 
discretization of an elliptic PDE on a two dimensional domain. 
The scheme decomposes the domain into thin subdomains, or ``slabs'' and uses a two-level approach
that is designed with parallelization in mind. 
The scheme takes advantage of $\mathcal H^2$-matrix structure emerging during factorization and
utilizes randomized algorithms to efficiently recover this structure.
As opposed to multi-level nested dissection schemes that incorporate the use of 
$\mathcal H$ or $\mathcal H^2$ matrices for a hierarchy of front sizes, 
SlabLU is a two-level scheme which only uses $\mathcal H^2$-matrix algebra for fronts of roughly
the same size. 
The simplicity allows the scheme to be easily tuned for
performance on modern architectures and GPUs.

\lsp 

The solver described is compatible with a range of different local discretizations, 
and numerical experiments demonstrate its performance for regular discretizations of rectangular and curved geometries.
The technique becomes particularly efficient when combined with very high-order accurate multi-domain spectral collocation schemes. 
With this discretization, a Helmholtz problem on a domain of size 
$1000 \lambda \times 1000 \lambda$ (for which $N=100 \rm{M}$) is solved 
in 15 minutes to 6 correct digits on a high-powered desktop with GPU acceleration.

\lsp

\end{minipage}

\end{center}

\section{Introduction}

\subsection{Problem setup}
We present a direct solver for boundary value problem of the form
\begin{equation}
\label{eq:bvp}
\left\{\begin{aligned}
\mathcal A u(x) =&\ f(x),\qquad&x \in \Omega,\\
u(x) =&\ g(x),\qquad&x \in \partial \Omega,
\end{aligned}\right.
\end{equation}
where $\mathcal A$ is a second order elliptic differential operator,
and $\Omega$ is a domain in two dimensions with boundary $\partial \Omega$.
The method works for a broad range of constant and variable coefficient differential
operators, but is particularly competitive for problems with highly oscillatory solutions
that are difficult to pre-condition.
For the sake of concreteness, we will focus on the case where $\mathcal A$ is a variable
coefficient Helmholtz operator
\begin{equation}
\label{eq:var_helm}
\mathcal A u(x) = -\Delta u(x) - \kappa^{2}b(x)u(x),
\end{equation}
where $\kappa$ is a reference (``typical'') wavenumber,
and where $b(x)$ is a smooth non-negative function.
Upon discretizing (\ref{eq:bvp}), one obtains a linear system
\begin{equation}
\label{eq:Au=f}
\mtx{A}\vct{u} = \vct{f},
\end{equation}
involving a coefficient matrix $\mtx{A}$ that is typically sparse. 
Our focus is on efficient algorithms for directly
building an invertible factorization of the matrix $\mtx{A}$. We specifically consider
two different discretization schemes, first a basic finite difference scheme with
second order convergence, and then a high (say $p=20$) order multidomain 
spectral 
collocation scheme \cite[Ch.~25]{2019_martinsson_book}. 
However, the techniques presented can easily be used with other (local)
discretization schemes such as finite element methods.

\subsection{Overview of proposed solver}
\label{sec:intro_overview}
The solver presented is based on a decomposition of the 
computational domain
into thin ``slabs'', as illustrated in Figure \ref{fig:domain}.
Unlike previously proposed sweeping schemes
\cite{2011_engquist_ying,gander2017restrictions,2013_sweep_demanet}
designed for preconditioning,
our objective is to directly factorize the coefficient matrix,
or at least compute a factorization that is sufficiently accurate
that it can handle problems involving strong backscattering.

To describe how the solver works, let us consider a simple model problem
where the PDE is discretized using a standard five-point finite difference
stencil on a uniform grid such as the one shown in Figure \ref{fig:domain}.
The nodes in the grid are arranged into slabs of width $b$, and are 
ordered as shown in Figure \ref{fig:domain}, resulting in a coefficient
matrix $\mtx{A}$ with the block diagonal sparsity pattern shown in 
Figure \ref{fig:sparseA}. The factorization of $\mtx{A}$ then proceeds through two stages.

\begin{figure}[!htb]
\captionsetup[subfigure]{labelfont=rm}
\centering
\begin{subfigure}[b]{0.45 \textwidth}
\centering
\tikzcuboid{16}{7}{0}{0.3}{4}
\caption{Original grid, partitioned into slabs
of width $b+2$, where $b = 3$.\\ \hspace{1em}}
\label{fig:domain}
\end{subfigure}%
\hfill
\begin{subfigure}[b]{0.45 \textwidth}
\centering
\reducedtikzcuboid{16}{7}{0}{0.3}{4}
\caption{Reduced grid, after eliminating blue nodes.
Only the red nodes are ``active''.\\}
\label{fig:reduced_grid}
\end{subfigure}

\begin{subfigure}[b]{0.45\textwidth}
    \centering
    \sparseA{18}{2.5}{4}{0.02}
    \caption{Sparsity structure of $\mtx A$ corresponding to the original grid. 
    Each block of $\mtx{A}$ is sparse.}
    \label{fig:sparseA}
\end{subfigure}%
\hfill
\begin{subfigure}[b]{0.45\textwidth}
    \centering
    \denseT{18}{2}{0.03}
    \vspace{2em}
    \caption{Sparsity structure of $\mtx T$ corresponding to the reduced grid. 
    Each block of $\mtx T$ is dense but has internal structure.}
    \label{fig:denseT}
\end{subfigure}
\caption{Illustration of the elimination order used in SlabLU.}
\end{figure}

In the first stage, the nodes that are internal to each slab 
(identified by the index vectors $I_2, I_4, \dots$ and shown as 
blue in Figure \ref{fig:domain}) are eliminated from the linear system, resulting in the reduced problem shown in
Figure \ref{fig:reduced_grid}, with the associated coefficient matrix shown
in Figure \ref{fig:denseT}. 
In this elimination step, we exploit that each subdomain is thin, which 
means that classical sparse direct solvers are particularly fast. 
To further accelerate this step, we use that the Schur complements that arise upon
the elimination of the interior nodes are rank structured. 
Specifically, they are ``HBS/HSS matrices'' \cite{2012_martinsson_FDS_survey,2013_efficient_FDS,2010_xia}
with \textit{exact} HBS/HSS rank at most $2b$. 
This allows us to accelerate this reduction step using a recently proposed 
randomized algorithm for 
compressing rank structured matrices \cite{levitt2024linear}.

The second stage is to factorize the remaining coefficient matrix $\mtx{T}$
shown in Figure \ref{fig:denseT}. This matrix is much smaller than the original
matrix $\mtx{A}$, but the sub-blocks are dense. Because we have efficiently formed
$\mtx T$ with $\mathcal H^2$-matrix structure, the reduced system can be
factorized in linear time for many elliptic PDEs (e.g. any coercive elliptic PDE, the steady state Stokes equation, and Helmholtz in the regime that the wavenumber is fixed as $N$ grows). 
With the use of $\mathcal H^2$-matrix algebra to factorize $\mtx T$, SlabLU requires linear time to 
store and factorize when the slab widths are chosen to be $\mathcal O(1)$.
In this work, we choose the slab width $b$ to grow slowly with the discretization size $N$
for performance considerations.

Our two-level framework offers a distinct advantage in terms of simplicity in both implementation and optimization. In the first stage,
we can leverage existing sparse direct solvers, which prove to be highly efficient, especially for thin 2D slabs. 
As we progress to the second stage, the fronts become larger in size, which may present a challenge when using traditional techniques.
The key benefit of SlabLU, in contrast to multi-level schemes, lies in the fact that we only need to develop specialized linear 
algebraic techniques for fronts of approximately the same size in the second stage. 
This stands in contrast to the necessity of developing such techniques for a hierarchy of front sizes in multi-level schemes.

We have found that for 2D problems, the dense 
operations are fast enough that exploiting rank structure to factorize $\mtx T$ is not
worthwhile when $N \leq 10^8$. 
Specifically, by choosing $b \sim \mathcal O(\sqrt{n})$, sparsity alone results in complexity $\mathcal O(N^{1.75})$ for the factorization stage, and $\mathcal O(N^{1.25})$ for the solve stage, when applied to 2D problems.
The simplicity of the two-level scheme allowed for parts of the factorization to 
be offloaded and accelerated on the GPU.
For meshes with 100 million points, the factorization can
be computed in 20 minutes on a desktop with an Intel i9-12900k CPU with 16 cores and an RTX 3090 GPU.
Once the factorization is available, subsequent solves take about a minute.
The numerical results feature timing results on a variety of architectures
to demonstrate that the scheme is portable to many hardware settings.

The scheme also interacts very well with high order discretization schemes such as those
described in \cite{2012_spectralcomposite} and \cite[Ch.~25]{2019_martinsson_book},
which makes it a particularly powerful tool for solving problems with highly oscillatory solutions.
The numerical results feature constant and variable coefficient 
Helmholtz problems on rectangular and curved domains. Using high order discretizations, we are able to 
discretize the PDE to 10 points per wavelength and accurately resolve
solutions on domains of size
$1000\lambda \times 1000\lambda$, where $\lambda$ is the wavelength,
to 6 digits of relative accuracy, compared to the true solution of the PDE.

\subsection{Context and related work}

Methods to solve (\ref{eq:Au=f}) can be characterized into two groups -- direct and iterative. 
The linear systems involved are typically ill-conditioned, which necessitates the use of specialized solvers. 
For problems with non-oscillatory solutions, multigrid methods are often highly effective
\cite{mccormick_multigrid,1987_ruge_algebraic,2017_xu_algebraic}.
For oscillatory problems, multigrid works less well \cite{ernst2012difficult}. 
Specialized preconditioners have been developed, and work well for many classes of problems, in particular those involving free-space problems \cite{gander2019class,gander2007optimized,erlangga2004class,erlangga2006comparison}.
In this context, the sweeping preconditioners of Engquist and Ying are of particular relevance, as they were an inspiration for the current work \cite{2011_engquist_ying,2011_engquist_ying_PML}.
However, oscillatory problems remain highly challenging to pre-condition, in particular in situations involving strong back-scattering, cavities, or problems trapped inside a finite domain.
Sparse direct methods, which factorize the matrix $\mtx A$ exactly, offer a robust means of solving challenging PDEs. 
They are also particularly advantageous in situations involving multiple right-hand sides or low-rank updates to the matrix $\mtx A$.

The solver we describe in this work is related to multi-frontal LU solvers \cite{2016_acta_sparse_direct_survey}
which often use a hierarchical nested dissection ordering of grid nodes \cite{amestoy1996approximate,george_1973}.
For a 2D grid with $N$ nodes,
the resulting techniques have complexity $\mathcal O(N^{3/2})$ to build and $\mathcal O(N\log N)$ complexity to solve, which is known
to be work optimal among solvers that exploit only sparsity in the system 
\cite{2006_davis_directsolverbook,1989_directbook_duff}. 
The $\mathcal H$ and $\mathcal H^2$-matrix algebras can be used to reduce the complexity of operations on dense matrices that arise in 
many contexts involving the discretization of integral equations and of PDEs \cite{2008_bebendorf_book,2003_borm_introduction_H_matrix,hackbusch2015hierarchical}. 
SlabLU is inspired by prior work in sparse direct solvers for PDEs which uses $\mathcal H$-matrices \cite{amestoy2017complexity,2018_chavez_accelerated,2016_ghysels_randomized_multifrontal,2009_xia_superfast,gillman2014direct,2017_pichon_sparse}. In particular, we were inspired by \cite{2018_chavez_accelerated}
which used a domain partitioning into planes (e.g. $b=1$) for a highly effective linear solver for low order discretizations.

A key feature of SlabLU is that unlike prior work, the rank structures that we exploit are \textit{exact}, relying only on the sparsity pattern of the original coefficient matrix (cf.~Section \ref{sec:rank_reduced}).
This makes the randomized compression particularly efficient, achieving very high computational efficiency with no loss of accuracy beyond floating point errors. 
Another novelty is the usage of a recently developed black box randomized algorithm for compressing rank structured matrices \cite{levitt2022linear} (which in turn draws on insights from \cite{2011_lin_lu_ying,2016_martinsson_hudson2,2011_martinsson_randomhudson}) when eliminating the interior nodes in each slab.

Importantly, the rank-deficiencies used for thin subdomains in SlabLU are present in both the non-oscillatory and oscillatory regimes.
Similar observations are used to develop
efficient solvers for integral equation discretizations on elongated domains \cite{1996_mich_elongated,2007_martinsson_elongated}.
In the general case, the interaction rank grows algebraically with the
wavenumber \cite{engquist2018approximate}, making the efficient use of $\mathcal H$ and $\mathcal H^2$ matrix techniques challenging, though it has been observed to work well in some situations \cite{banjai2008hierarchical,betcke2017computationally,2011_xia_dehoop}. For the purposes of this work which focuses on 2D domains, 
we only use $\mathcal H^2$-matrix compression to form the reduced system $\mtx T$ efficiently,
then use highly efficient dense linear algebra routines to factorize the reduced system, an approach which we demonstrate to be effective
in the numerical results section.

%To resolve the oscillatory solutions of the Helmholtz equation, the equation must be discretized to a constant number of points per wavelength, roughly 10. PDE formulations suffer from the so-called ``pollution effect'' \cite{1997_babuska_pollution} --- as the wavenumber increases, one may suffer from phase errors unless you discretize with \textit{even more} points per wavelength. To combat these effects, we use a multi-domain high-order spectral collocation scheme, called the ``Hierarchical Poincar\'e-Steklov (HPS)'' scheme, described in \cite{babb2018accelerated,2019_martinsson_book,hao2016direct}.

%\begin{remark}
%Fast direct solvers for the Helmholtz equation that attain acceleration by exploiting 
%rank structure \pgmnote{Citations.} generally 
%lose efficiency in the high frequency regime, as the numerical 
%ranks grow with the wavenumber. 
%In contrast, the rank structures that we exploit in the first stage rely on exact algebraic
%properties, which means that they work for all wavenumbers.
%\end{remark}

\subsection{Extensions and limitations}

The solver presented is purely algebraic and can be applied to a range of different discretization schemes, including finite element and finite volume methods. 
In this manuscript, we restrict attention to regular discretizations of domains that are either rectangular themselves, 
or can be mapped smoothly to a union of elongated rectangles or slabs.
It is possible to adapt the method to more general discretizations with local refinement,
so long as it is simple to partition the computational domain into index sets corresponding to elongated slabs.
For some high order discretizations (e.g. high order finite differences), widening stencils may lead to large pre-factors when using
SlabLU, though this is also a challenge for sparse direct solvers in general \cite[Ch.~20]{2019_martinsson_book}.

While we in this manuscript restrict attention to the two dimensional case, the method is designed to handle three dimensional problems as well. 
All ideas presented carry over directly, but additional complications do arise. 
The key challenge is that in three dimensions, it is no longer feasible to use dense linear algebra when factorizing the block tridiagonal reduced coefficient matrix $\mtx{T}$.
However, the 3D version of SlabLU is also very easy to parallelize, and the idea of using randomized compression combined with efficient sparse direct solvers to eliminate the nodes interior to each slab still applies, cf.~Section \ref{sec:conc}.

\section{Discretization and node ordering}

We introduce two different discretization techniques for
(\ref{eq:bvp}). The first is simply the standard second order accurate
five point finite difference stencil. Since this discretization is very 
well known, it allows us to describe how the solver works without the need
to introduce cumbersome background material. To demonstrate that the solver
works for a broader class of discretization schemes, the numerical experiments
reported in Section \ref{sec:num} also include results that rely on the high 
order accurate \textit{Hierarchical Poincar\'e-Steklov (HPS)} scheme, which
we briefly describe in Section \ref{sec:other_disc}.

\subsection{A model problem based on the five point stencil}
\label{sec:modeldisc}

For purposes of describing the factorization scheme, let us introduce a very
simple discretization of the boundary value problem (\ref{eq:bvp}).
We work 
with a rectangular domain $\Omega=[0,L_1] \times [0,L_2]$ and the second order
linear elliptic operator $\mathcal A$ defined by (\ref{eq:var_helm}). We assume 
that $L_1 \ge L_2$, and that $L_{1} = hn_{1}$ and $L_{2} = hn_{2}$ for some grid 
spacing $h$ and some positive integers $n_{1}$ and $n_{2}$.
We then
discretize $\mathcal{A}$ with a standard second-order finite difference scheme,
to obtain the linear system
\begin{equation}
\label{eq:discretiation}
\frac 1 {h^2} 
\bigl(
\vct{u}(n_{\rm w}) + 
\vct{u}(n_{\rm e}) +
\vct{u}(n_{\rm n}) +
\vct{u}(n_{\rm s}) -
4 \vct{u}(n)
\bigr)
- \kappa^2 \vct{b}(n)\vct{u}(n)
=
\vct{f}(n).
\end{equation}
The vector $\vct{f}$ holds values of the body load at the discretization nodes, and the vector $u$ holds approximations to the solution $u$. 
See Figure \ref{fig:fivept_stencil} for a visualization of the 5 point stencil. 
We write the system (\ref{eq:discretiation}) compactly as
$\mtx{A}\vct{u} = \vct{f}$.

\subsection{Clustering of the nodes}
\label{sec:clustering}

\begin{wrapfigure}{r}{0.30\textwidth}
\vspace{0pt}
\centering
% No need for \captionsetup{type=figure} inside wrapfigure; just use \caption.
% Make label text the same size as the main body:
\tikzset{
  every label/.style={font=\normalsize, inner sep=2pt, label distance=1.5pt},
  smallpoint/.style={circle, fill=black, draw=black, inner sep=0.7pt, minimum size=4pt}
}

% Node command: accepts style, coordinates, node name, label text, and label position
\newcommand{\stencilpt}[5][]{\node[smallpoint, #1, label={#5:#4}] at (#2) (#3) {};}

% Scale the geometry only; text stays at body size.
\begin{tikzpicture}[scale=0.9, baseline=(current bounding box.north)]
    \stencilpt{ 0,0}{n}{$\vct {u}(n)$}{below right};
    \stencilpt{ 0,1.25}{n_n}{$\vct {u}(n_{\rm n})$}{above};
    \stencilpt{ 0,-1.25}{n_s}{$\vct {u}(n_{\rm s})$}{below};
    \stencilpt{ 1.25,0}{n_e}{$\vct {u}(n_{\rm e})$}{above};
    \stencilpt{ -1.25,0}{n_w}{$\vct {u}(n_{\rm w})$}{above};
    \draw (n) -- (n_n) (n) -- (n_e) (n) -- (n_s) (n) -- (n_w);
\end{tikzpicture}

\caption{\centering Five-point \newline stencil in 2D.}

\label{fig:fivept_stencil}
\end{wrapfigure}

We next subdivide the computational domain into thin ``slabs'', as shown in 
Figure \ref{fig:domain}. We let $b$ denote the number of grid points in each
slab ($b=4$ in Figure \ref{fig:domain}), and then introduce index vectors $I_{1},\, I_{2},\, I_{3}, \dots$ that keep
track of which slabs each grid point belongs to. The odd numbered index
vectors $I_{1},\, I_{3},\,I_{5},\dots$ indicate nodes on the interfaces between
slabs (red in Figure \ref{fig:domain}), while the even numbered ones indicate
nodes that are interior to each slab (blue in Figure \ref{fig:domain}). 
With this ordering of the grid points, the stiffness matrix associated with the
discretization (\ref{eq:discretiation}) has the sparsity pattern shown in 
Figure \ref{fig:sparseA}.

\subsection{High order discretizations}
\label{sec:other_disc}
To accurately resolve oscillatory wave phenomena, we rely on a high order accurate
multi-domain spectral collocation discretization known as the Hierarchical 
Poincar\'e-Steklov scheme (HPS).
This discretization scheme is designed to allow for high choices of the local discretization order $p$ without degrading the performance of direct solvers.
In HPS, the computational domain is subdivided into small subdomain, and a $p\times p$ tensor product grid of Chebyshev nodes is placed on each subdomain.

\begin{figure}[!htb]
\fbox{
\centering
\begin{subfigure}{0.40\textwidth}
    \centering
    \includegraphics[width=0.8\textwidth]{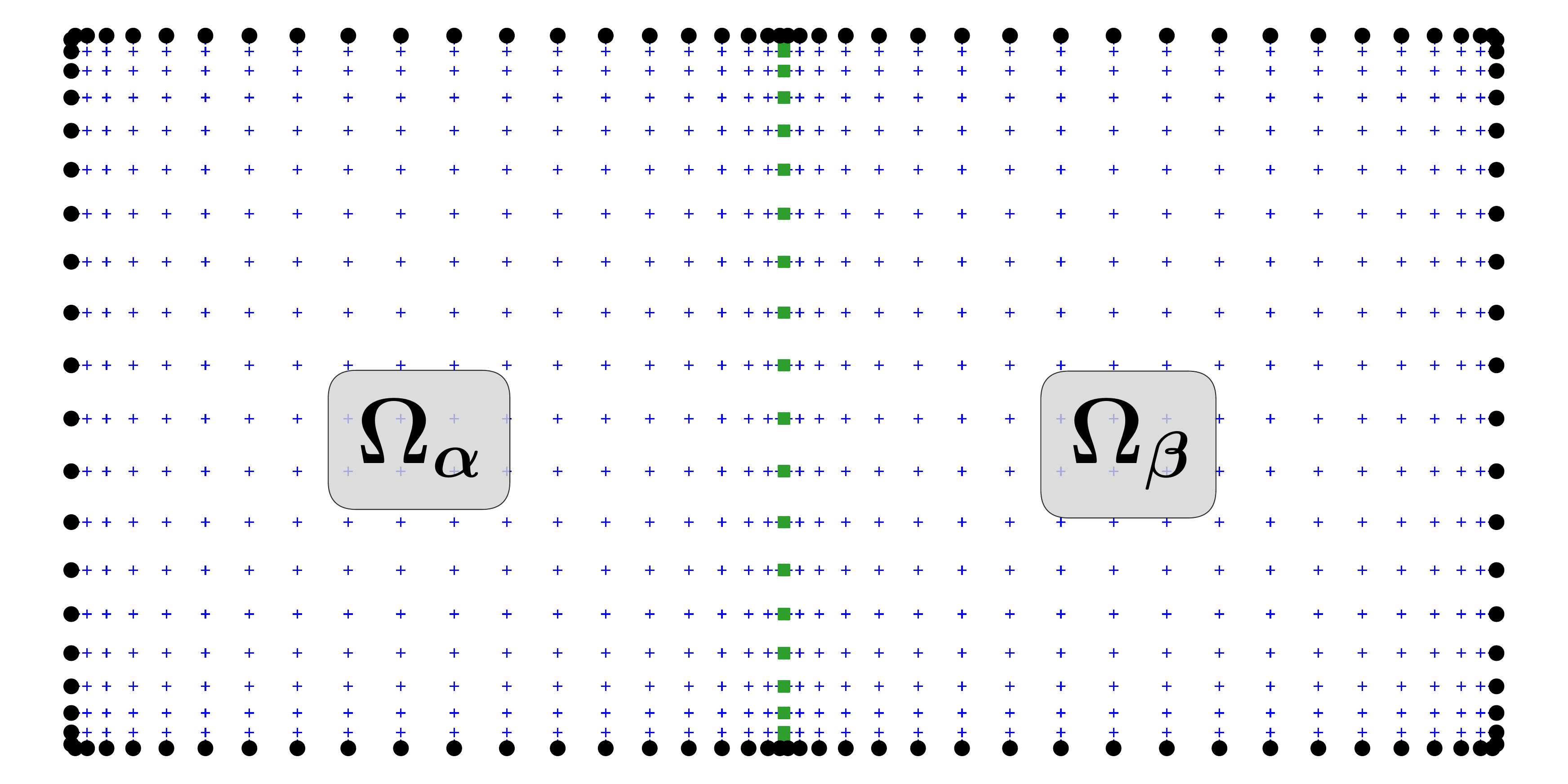}
    Consider the discretization of $\Omega = \Omega_\alpha \cup \Omega_\beta$
    with labeling of the nodes $I = I_1 \cup I_2 \cup I_3$ so that 
    \begin{align*}
    I_1:\ &\text{collocation of (\ref{eq:Au=f}) on}\ \Omega_\alpha\\
    I_2:\ &\text{collocation of (\ref{eq:Au=f}) on}\ \Omega_\beta\\
    I_3:\ &\text{shared boundary, where }\\
    &\text{continuity of $\partial u / \partial x$ enforced}
    \end{align*}
\end{subfigure}
\hfill
\begin{subfigure}{0.58\textwidth}
    The stiffness matrix can be partially factorized to decouple $I_1$ and $I_2$ from the rest of the system as
    \[
    \left[
    \begin{array}{@{}ccc|ccc|c@{}}
    &&&&&&\\
    & \mtx A_{11} & & & & &\mtx A_{13}\\
    &&&&&&\\ \hline
    &&&&&&\\
    & & & & \mtx A_{22} & & \mtx A_{23}\\
    &&&&&&\\ \hline
    & \mtx D^{\alpha}_{31} &&& -\mtx D^{\beta}_{32} && \mtx N
    \end{array}
    \right]
    =
    \renewcommand{\arraystretch}{1.18}
    \mtx L
    \left[
    \begin{array}{@{}c|c|c@{}}
     \mtx A_{11} & \\
    \hline
    & \mtx A_{22} & \\ \hline
    && \mtx {\tilde A}_{33}
    \end{array}
    \right]
    \mtx U,
    \renewcommand{\arraystretch}{1.0}
    \]
    where 
    \begin{align*}
    \mtx N &= \mtx D^{\alpha}_{33}-\mtx D^{\beta}_{33}\ \text{enforces continuity of $\partial u/ \partial x$},\\
    \mtx {\tilde A}_{33} &= \mtx N - \mtx D^{\alpha}_{31} \mtx A_{11}^{-1} \mtx A_{13} + \mtx D^{\beta}_{32} \mtx A_{22}^{-1} \mtx A_{23}.
    \end{align*}
\end{subfigure}
}
\caption{The figure above provides an illustration of static condensation for a simple HPS discretization of two subdomains for (\ref{eq:Au=f}) with
Dirichlet data prescribed. See \cite[~Ch. 24]{2019_martinsson_book} for further details on the details of the discretization. Static condensation is the process of partially factorizing $\mtx A$ and decoupling the nodes internal to each
subdomain from the rest of the system. The reduced system (\ref{eq:reduced_hps_system}) is on the interfaces between subdomains.
}
\label{fig:hps_two_sub}
\end{figure}

As a brief illustration of the discretization scheme, consider a partitioning of $\Omega$ into two subdomain, cf. Figure \ref{fig:hps_two_sub}.
For the nodes internal to each subdomain, we discretize (\ref{eq:Au=f}) through collocation of the spectral differentiation operator.
For the nodes on boundaries between subdomains, we enforce continuity of the normal derivatives. 

\begin{figure}[!htb]
\centering
\begin{minipage}{0.4\textwidth}
\centering
\includegraphics[width=\textwidth]{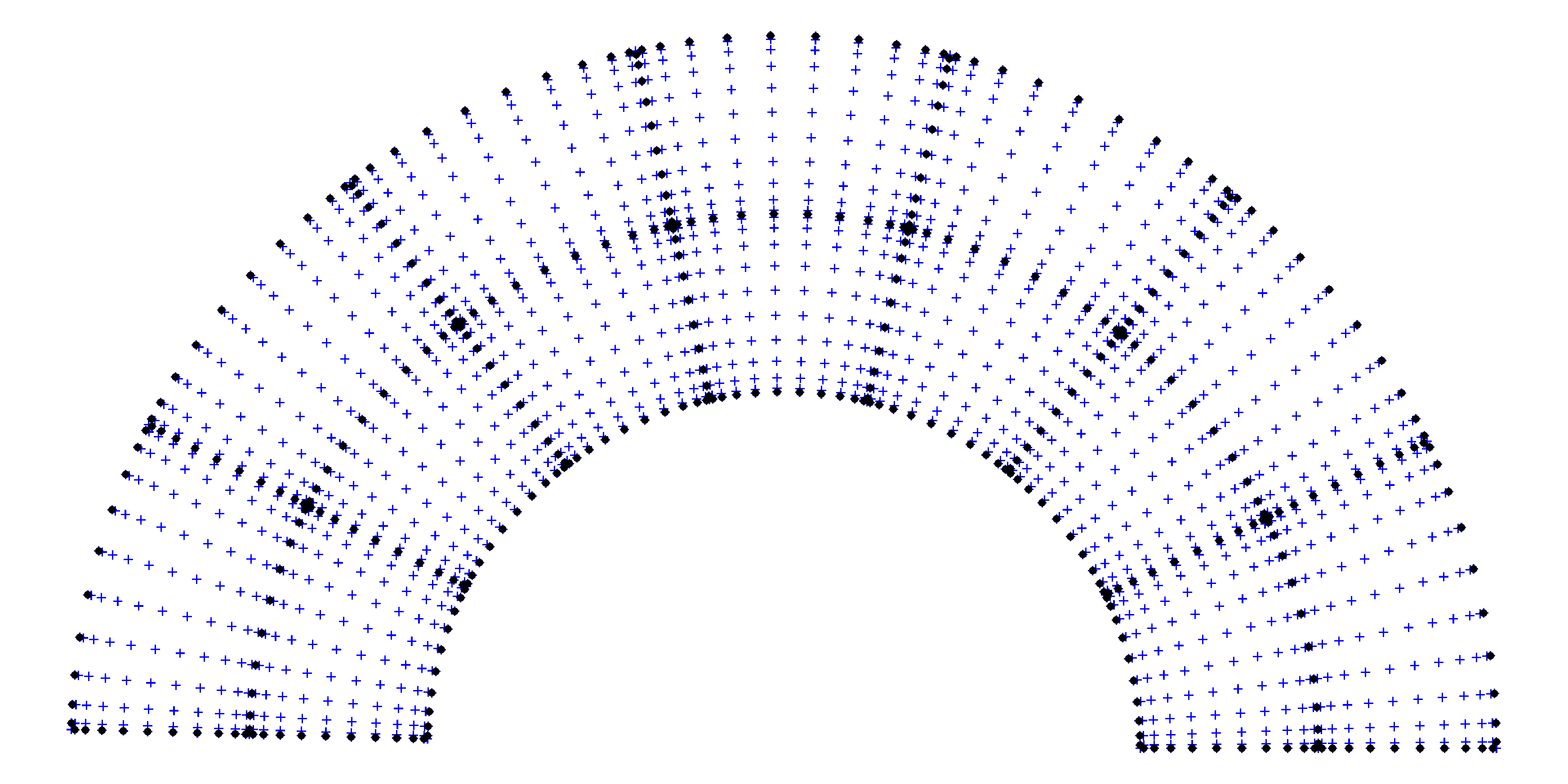}
\end{minipage}%
$\Rightarrow$
\begin{minipage}{0.4\textwidth}
\centering
\includegraphics[width=\textwidth]{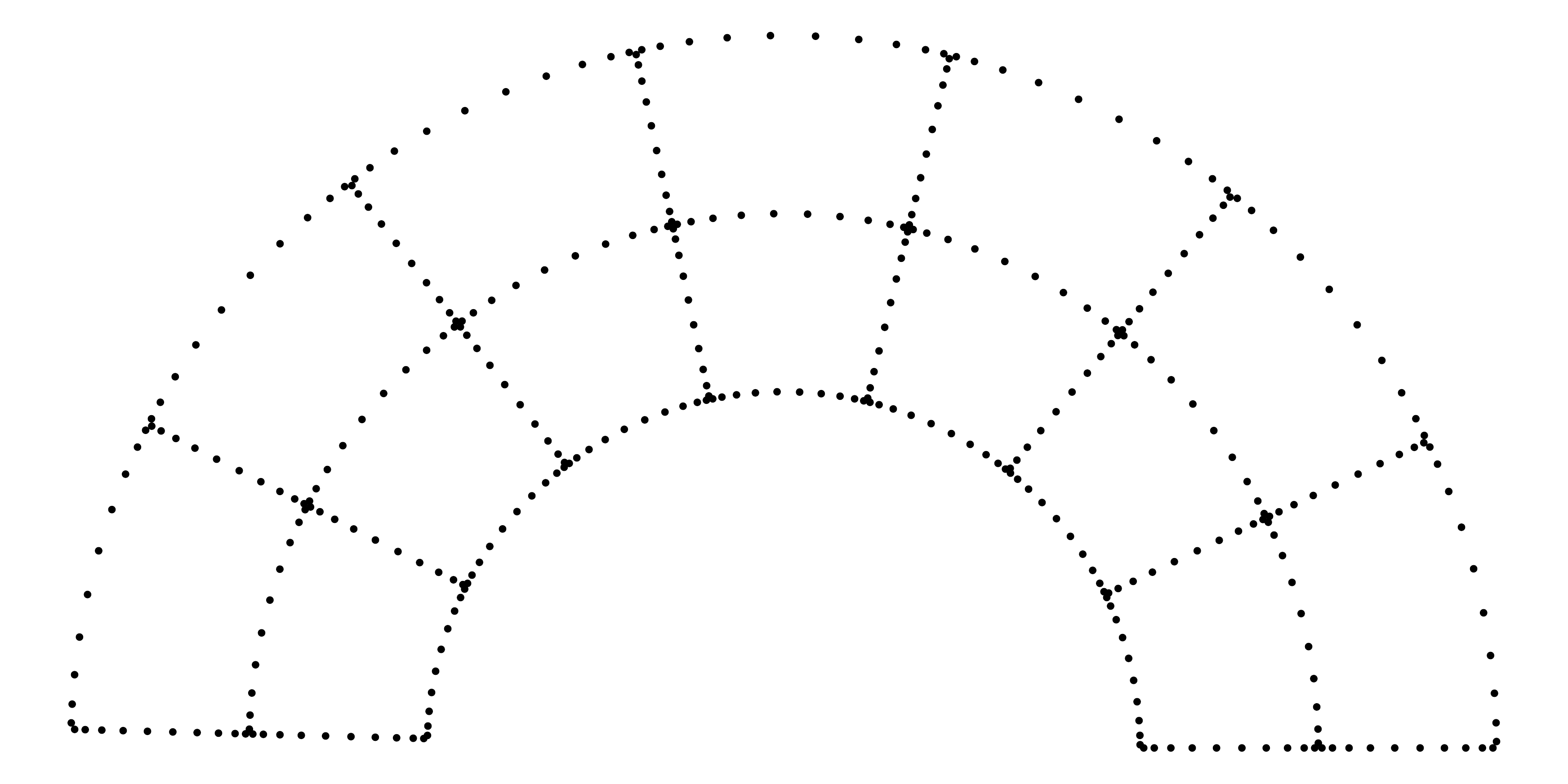}
\end{minipage}
\caption{\small HPS is a multi-domain spectral collocation scheme where the PDE is enforced on each subdomain interior using dense spectral differentiation. Prior to interfacing with SlabLU, we ``eliminate'' the interior blue nodes in parallel and produce an equivalent system to solve on the boundaries. The original grid has $n_1 \times n_2$ points, and remaining grid has $\approx n_1 n_2/p$ points.}
\label{fig:hps_reduced_grid}
\end{figure}

To improve efficiency when HPS is combined with sparse direct solvers, we ``eliminate'' the nodes internal to each subdomain 
by partially computing an LU decomposition for the internal nodes that decouples them from the rest of system. This process 
is known as static condensation and leads a reduced system
\begin{equation}\
\mtx {\tilde A} \mtx {\tilde u} = \mtx {\tilde f}
\label{eq:reduced_hps_system}
\end{equation}
with modified interactions between subdomain interfaces. 
The discretization can be generalized to any domain that can be smoothly mapped to a union of square subdomains, cf. Figure \ref{fig:hps_reduced_grid}.
For further details, see \cite[Ch.~25]{2019_martinsson_book}, 
as well as \cite{babb2018accelerated,2013_martinsson_ItI,hao2016direct,2012_spectralcomposite,beams2020parallel,2019_gillman_HPS_adaptive}.

Because spectral differentiation is a dense operator, static condensation requires $\mathcal O(p^4 N)$ flops and the solution operators for each subdomain 
require $\mathcal O (p^2 N)$ bits to store. In our implementation, we use the fact that the leaf operations can be done in parallel
and attain high performance on the GPU with batched linear algebra, which we demonstrate to be an effective approach for $p$ up to 42 in \cite{yesypenko2022parallel}. 
The leaf operations are so efficient that we discard the factorized operators on leaf nodes to save on space and refactorize as needed, cf. Section \ref{sec:complexity_hps}.

For efficiently handling higher orders of $p$, one can use ultraspherical polynomials \cite{2013_olver_townsend_wellcond} instead of Chebyshev polynomials 
which lead to sparse and well-conditioned differentiation on the leaf nodes, see \cite{2020_spectral_townsend_fortunato,2021_fortunato_ultraspherical,aurentz2020symmetrizing}
for efficient multi-domain spectral methods on 2D volumes and \cite{fortunato2022_surface} for methods on surface PDEs.

\begin{remark}
A key point of the present work is that the solver has only \textit{two} levels,
which makes the ``H'' in ``HPS'' a slight misnomer, as it refers to ``hierarchical''. 
We nevertheless stick with the ``HPS'' acronym to conform with the prior literature.
\end{remark}

\section{Stage One: Elimination of nodes interior to each slab}
\label{sec:stageone}

This section describes the process that we use to eliminate 
the nodes interior to each slab that we sketched out in Section \ref{sec:intro_overview}.
The objective is to reduce the sparse stiffness matrix $\mtx{A}$ 
(illustrated in Figure \ref{fig:sparseA}) into the smaller block
tridiagonal matrix $\mtx{T}$ (illustrated in Figure \ref{fig:denseT}).
The techniques described form the core algorithmic innovation
of the manuscript.

\subsection{Schur complements}
With the ordering introduced in Section \ref{sec:clustering}, the
coefficient matrix $\mtx{A}$ has the block form
\begin{equation}
\label{eq:dell1}
\begin{bmatrix}
\mtx A_{11} & \mtx A_{12} & \mtx 0 & \mtx 0 & \mtx 0 & \dots \\
\mtx A_{21} & \mtx A_{22} & \mtx A_{23} & \mtx 0 &  \mtx 0 & \dots \\
\mtx 0 & \mtx A_{32} & \mtx A_{33} & \mtx A_{34} &   \mtx 0 & \dots\\
\mtx{0} & \mtx 0 & \mtx A_{43} & \mtx A_{44} & \mtx A_{45} & \dots\\
\vdots & \vdots & \vdots & \vdots & \vdots & \vdots
\end{bmatrix}
\begin{bmatrix} \mtx u_1\\\mtx u_2\\ \mtx u_3 \\ \mtx u_4 \\ \vdots \end{bmatrix} =
\begin{bmatrix} \mtx f_1\\ \mtx f_2\\ \mtx f_3\\ \mtx f_4 \\ \vdots        \end{bmatrix}.
\end{equation}
We eliminate the vectors $\vct{u}_{2},\,\vct{u}_{4},\,\vct{u}_{6},\,\dots$
that represent unknown variables in the interior of each slab through a 
step of block Gaussian elimination. To be precise, we insert the relation
\begin{equation}
\vct{u}_{i} = 
\mtx{A}_{ii}^{-1}\bigl(\vct{f}_{i} - \mtx{A}_{i,i-1}\vct{u}_{i-1} - \mtx{A}_{i,i+1}\vct{u}_{i+1}\bigr),
\qquad i = 2,\,4,\,6,\,\dots
\label{eq:usol_even}
\end{equation}
into the odd-numbered rows in (\ref{eq:dell1}) to obtain the reduced system
\begin{equation}
\label{eq:dell3}
\begin{bmatrix}
\mtx T_{11} & \mtx T_{13} & \mtx 0 & \mtx 0 & \mtx 0 & \dots \\
\mtx T_{31} & \mtx T_{33} & \mtx T_{35} & \mtx 0 &  \mtx 0 & \dots \\
\mtx 0 & \mtx T_{53} & \mtx T_{55} & \mtx T_{57} &   \mtx 0 & \dots\\
\mtx{0} & \mtx 0 & \mtx T_{57} & \mtx T_{77} & \mtx T_{79} & \dots\\
\vdots & \vdots & \vdots & \vdots & \vdots & \vdots
\end{bmatrix}
\begin{bmatrix} \mtx u_1\\\mtx u_3\\ \mtx u_5 \\ \mtx u_7 \\ \vdots \end{bmatrix} =
\begin{bmatrix} \tilde{\vct{f}}_{1} \\ \tilde{\vct{f}}_{3} \\\tilde{\vct{f}}_{5} \\ \tilde{\mtx f}_7\\ \vdots \end{bmatrix},
\end{equation}
where the sub-blocks of $\mtx T$ are defined as
\begin{align}
\label{eq:T11}
\mtx{T}_{11} &= \mtx A_{11} - \mtx{A}_{ 12}\ \mtx{A}_{ 22}^{-1}\ \mtx{A}_{ 21}, \\
\label{eq:T13}
\mtx{T}_{13} &= \mtx A_{13} - \mtx{A}_{ 12}\ \mtx{A}_{ 22}^{-1}\ \mtx{A}_{ 23}, \\
\label{eq:T31}
\mtx{T}_{31} &= \mtx A_{31} - \mtx{A}_{ 32}\ \mtx{A}_{ 22}^{-1}\ \mtx{A}_{ 23}, \\
\label{eq:T33}
\mtx{T}_{33} &= \mtx A_{33} - \mtx{A}_{ 32}\ \mtx{A}_{ 22}^{-1}\ \mtx{A}_{ 23}
                              - \mtx{A}_{ 34}\ \mtx{A}_{ 44}^{-1}\ \mtx{A}_{ 43}, \\
\label{eq:T35}
\mtx{T}_{35} &= \mtx A_{35} - \mtx{A}_{ 34}\ \mtx{A}_{ 44}^{-1}\ \mtx{A}_{ 23},
\end{align}
and so on. The reduced right-hand sides $\mtx{\tilde f}$ are defined as
\begin{align}
\label{eq:f1}
\mtx{\tilde f}_{ 1} &= \mtx f_{ 1} - \mtx{A}_{ 12}\ \mtx{A}_{22}^{-1}\ \mtx f_2, \\
\label{eq:f3}
\mtx{\tilde f}_{ 3} &= \mtx f_{ 3} - \mtx{A}_{ 32}\ \mtx{A}_{22}^{-1}\ \mtx f_2 
                                   - \mtx{A}_{ 34}\ \mtx{A}_{44}^{-1}\ \mtx f_4, \\
\label{eq:f5}
\mtx{\tilde f}_{ 5} &= \mtx f_{ 5} - \mtx{A}_{ 54}\ \mtx{A}_{44}^{-1}\ \mtx f_4
                                   - \mtx{A}_{ 56}\ \mtx{A}_{66}^{-1}\ \mtx f_6,
\end{align}
etc.

\subsection{Rank structure in the reduced blocks}
\label{sec:rank_reduced}

We next discuss algebraic properties of the blocks of the reduced coefficient matrix that allow for $\mtx T$ to be formed efficiently. The sub-blocks of $\mtx T$ are Schur-complements of sparse matrices. For instance, block $\mtx T_{\rm 11}$
has the formula
\begin{equation}\underset{n_2 \times n_2}{\mtx T_{11}} = \underset{n_2 \times n_2}{\mtx A_{11}} - \underset{n_2 \times n_2 b} {\mtx A_{12}}\ 
\underset{n_2 b \times n_2 b} {\mtx A_{22}^{-1}}\ \underset{n_2 b \times n_2} {\mtx A_{21}},
\label{eq:T11_detailed}
\end{equation}
where $\mtx A_{11}, \mtx A_{12}, \mtx A_{21}$ are sparse with $\mathcal O(n_2)$ non-zero entries and $\mtx A_{22}$ is a sparse banded matrix. 
Factorizing $\mtx A_{22}$ can be done efficiently using a sparse direct solver, 
but forming $\mtx A_{22}^{-1} \mtx A_{21}$ naively may be slow and memory-intensive, especially considering that 
$\mtx A_{21}$ is sparse and would need to be converted to dense to interface with solve routines. 

We use an alternate approach for efficiently forming $\mtx T_{11}$ that achieves high arithmetic intensity
while maintaining a very low memory footprint. First, we prove algebraic properties about $\mtx T_{11}$, which is a
dense but structured matrix that only needs $\mathcal O(n_2 b)$ storage 
in exact precision. 
The matrix $\mtx T_{11}$ is compressible in a format called Hierarchically Block-Separable (HBS) or Hierarchically
Semi-Separable (HSS) with exact HBS rank at most $2b$. 
HBS matrices are a type of hierarchical matrix ($\mathcal H^2$-matrix), which allow dense matrices to be stored efficiently
by exploiting low-rank structure in sub-blocks at different levels of granularity \cite{2008_bebendorf_book,hackbusch2015hierarchical,2019_martinsson_book}.
The sub-blocks of $\mtx T$ are also compressible in $\mathcal H$-matrix formats, e.g. Hierarchical Off-Diagonal Low Rank (HODLR). 
The rank property of $\mtx T_{11}$ is formally stated in Proposition \ref{prop:rank_property}.

After establishing the HBS structure of $\mtx T_{11}$, we then 
describe how this structure can be recovered with only $\mathcal O(b)$ matrix vector products
of $\mtx T_{11}$ and $\mtx T_{11}^*$. Vectors can efficiently be applied 
because $\mtx T_{11}$ and its transpose are compositions of sparse matrices.
Since $\mtx T_{11}$ is admissible as a HODLR matrix, the structure can be efficiently recovered from matrix-vector products in this format as well
 \cite{2011_lin_lu_ying,2016_martinsson_hudson2}, though more vectors are required for reconstruction in this format.

\begin{figure}
\begin{subfigure}{0.4\textwidth}
\centering
\simpleproofpicture{30}{5}{0.15}
\subcaption{Contiguous set of points $J_B \subset I_1$ as described in Proposition \ref{prop:rank_property}.}
\label{fig:box_farfield}
\end{subfigure}
\hfill
\begin{subfigure}{0.58\textwidth}
\centering

\vspace{2.5em}

\begin{tikzpicture}[scale=0.17]

\pgfmathsetmacro{\b}{2}
\pgfmathsetmacro{\k}{1}

\pgfmathsetmacro{\m}{6}

\pgfmathsetmacro{\offsetx}{2.5*\m}
\pgfmathsetmacro{\offsety}{0}

% low rank row
\filldraw[fill=gray, thick] (\offsetx,\b) rectangle (\offsetx+\m*\b,\m*\b+\b);

% horizontal lines separating rows
\foreach \x in {1,...,\m}
\draw[thick, white] (\offsetx,\offsety+\x*\b) -- (\offsetx+2*\m,\offsety+\x*\b); 
\filldraw[fill=none, thick] (\offsetx,\offsety+\b) rectangle (\offsetx+\m*\b,\offsety+\m*\b+\b);

\foreach \x in {1,...,\m}
\draw[thick,white,fill=darkgray] (\offsetx + \m*\b - \x * \b + \b, \offsety + \x * \b + \b) rectangle (\offsetx + \m*\b - \x * \b , \offsety + \x * \b );

\filldraw[fill=none, thick] (\offsetx +0,\offsety+\b) rectangle (\offsetx +\m*\b,\offsety+\m*\b+\b);

% low rank column
\pgfmathsetmacro{\offsetx}{5*\m}
\pgfmathsetmacro{\offsety}{0}

\filldraw[fill=gray, thick] (\offsetx +0,\offsety+\b) rectangle (\offsetx +\m*\b,\offsety+\m*\b+\b);

% vertical lines separating columns
\foreach \x in {1,...,\m}
\draw[thick, white] (\offsetx + \x *\b,\offsety+\b) -- (\offsetx + \x * \b,\offsety+2*\m+\b);

\foreach \x in {1,...,\m}
\draw[thick, white,fill=darkgray] (\offsetx + \m*\b - \x * \b + \b, \offsety+\x * \b + \b) rectangle (\offsetx+\m*\b - \x * \b , \offsety+\x * \b );

\filldraw[fill=none, thick] (\offsetx +0,\offsety+\b) rectangle (\offsetx +\m*\b,\offsety+\m*\b+\b);

\pgfmathsetmacro{\offsetx}{7.5*\m}
\pgfmathsetmacro{\offsety}{0}

\draw[fill=darkgray] (\offsetx,\m+1.5*\b) rectangle (\offsetx+\b,\m+2.5*\b);
\draw[fill=gray] (\offsetx,\m ) rectangle (\offsetx+\b, \m+ 1.0*\b);
\node at (\offsetx+4*\b,\m+2*\b) {full rank};
\node at (\offsetx+4*\b,\m+0.5*\b) {low rank};

\end{tikzpicture}

\vspace{4em}

\subcaption{The subblock $\mtx T_{11}$ is compressible as Hierarchically Block-Separable matrices (HBS) and can be 
tesellated so that sub-blocks are low-rank or small enough to be stored densely.}
\label{fig:hbs}
\end{subfigure}
\caption{The geometry of slab interface $I_1$ is shown in Figure \ref{fig:box_farfield}. The submatrix $\mtx T_{11}$ is compressible
as an HBS matrix (cf. Figure \ref{fig:hbs} for an illustration).}
\label{fig:hbs_slab_twofigs}
\end{figure}

\vspace{0.5em}

\begin{restatable}[Rank Property]{proposition}{rankprop}
\label{prop:rank_property}
Let $J_B$ be a contiguous set of points on the slab interface $I_1$, and 
let $J_F$ be the rest of the points $J_{F} = I_1 \setminus J_{B}$.
The sub-matrices
$(\mtx T_{11})_{BF},\ (\mtx T_{11})_{FB}$ of the matrix $\mtx T$ have exact rank at most $2b$. 
\end{restatable}

\vspace{0.5em}

See Figure \ref{fig:hbs_slab_twofigs} for an illustration of the slab interface and the resulting structure the sub-blocks of $\mtx T$, which are defined in (\ref{eq:T11}-\ref{eq:T35}).
The proof is in Appendix \ref{sec:appendix}.

\subsection{Recovering {$\mathcal{H}^2$}-matrix structure from matrix-vector products}
\label{sec:hmatrix_rand}
We next describe how to extract an $\mathcal H^2$-matrix representation 
of the reduced blocks purely from matrix-vector products.
For concreteness, we are trying to recover $\mtx T_{11} \in \mathbb{R}^{n_2 \times n_2}$ 
as an HBS matrix with HBS rank at most $2b$. Typically, $\mathcal H^2$-matrices are used 
when the cost of forming or factorizing these matrices densely, is prohibitively large. 
In the context of SlabLU, $\mtx T_{11} \in \mathbb{R}^{n_2 \times n_2}$ can be stored densely for the problem sizes of interest, but 
traditional methods for forming $\mtx T_{11}$ densely may be inefficient, 
as described in Section \ref{sec:rank_reduced}. 
Instead, we recover $\mtx T_{11}$ as an HBS matrix with HBS rank at most $2b$ from matrix-matrix products 
\begin{equation}\underset{n_2 \times s}{\mtx Y} = 
\underset{n_2 \times n_2}{\mtx T_{11}}\ \underset{n_2 \times s}{\mtx \Omega}, 
\qquad \underset{n_2 \times s}{\mtx Z} = \underset{n_2 \times n_2}{\mtx T_{11}^*}\ 
\underset{n_2 \times s}{\mtx \Psi},
\label{eq:rand_samples}
\end{equation} 
where $\mtx \Omega, \mtx \Psi$ are Gaussian random matrices and $s = 6b$ using the 
algorithm presented in \cite{levitt2022linear}. This is theoretically possible 
because an HBS matrix of size $n_2 \times n_2$ with HBS rank at most $2b$ can be encoded in $\mathcal O(n_2 b)$ storage.
The HBS structure can be recovered from samples $\mtx Y, \mtx Z$ after 
post-processing in $\mathcal O(n_2 b^2)$ flops. 
The algorithm presented in \cite{levitt2022linear} can be seen as an extension of algorithms for 
recovering low-rank factors from random sketches \cite{martinsson2020randomized}. 
A particular advantage of these algorithms is that they scale linearly and are truly black-box.
The matrix-matrix products (\ref{eq:rand_samples}) are simple to evaluate using the matrix-free 
formula (\ref{eq:T11_detailed}) of $\mtx T_{11}$ and its transpose. 
Applying $\mtx T_{11}$ involves
two applications of sparse matrices and two triangular solves using pre-computed
sparse triangular factors.

\section{Stage Two: Factorizing the reduced block tridiagonal coefficient matrix}
\label{sec:stagetwo}

The elimination of nodes interior to each slab that we described in 
Section \ref{sec:stageone} results in a reduced linear system
\begin{equation}
{\mtx T} \mtx {\tilde u} = \mtx {\tilde f}, 
\label{eq:Tred_system}
\end{equation}
where $\mtx {\tilde u}$ is the reduced solution vector and where $\mtx {\tilde f}$ is the equivalent body load on slab interfaces $I_1, I_3, \dots$ 
The elimination process described in Section 
\ref{sec:stageone} results a block-tridiagonal reduced system (\ref{eq:Tred_system}) with sub-blocks 
that are compressible in HBS format.
Solving a system involving a block tridiagonal matrix is straight-forward
using a blocked version of Gaussian elimination, described in Algorithms \ref{alg:sweeping} and \ref{alg:sweeping_solve}. 

\RestyleAlgo{boxruled}
\begin{algorithm}[!htb]
  \caption{Sweeping build stage\label{alg:sweeping}.\\
  Given a block-tridiagonal matrix $\mtx T$, Algorithm \ref{alg:sweeping} builds a direct solver for $\mtx T$, with the result stored in the matrices $\mtx S_1, \mtx S_3, \dots$}
  $\mtx S_1 \gets \mtx T_{11}$\;
  Compute and store $\mtx S_1^{-1}$\;
  \For{$j=2,\dots,n_1/b$}{
  $\mtx S_{ 2j+1} \gets \mtx T_{ 2j+1,2j+1} - \mtx T_{ 2j+1,2j-1} \mtx S_{2j-1}^{-1} \mtx T_{ 2j-1,2j+1}$\;
  Compute and store $\mtx S_{ 2j+ 1}^{-1}$\;
  }
  \end{algorithm}
  \RestyleAlgo{boxruled}
\begin{algorithm}[!htb]
  \caption{Sweeping solve.\label{alg:sweeping_solve}\\
  Given a body load $\mtx {\tilde f}$ and precomputed inverses of $\mtx S_1, \mtx S_3, \dots$,
  Algorithm \ref{alg:sweeping_solve} computes the solution vector $\mtx {\tilde u}$ by overwriting the original vector in place.}
  \For{$j=1,\dots, n_1/b$}{
  $\mtx {\tilde f}_{ 2j + 1} \gets \mtx {\tilde f}_{ 2j+1} - \mtx T_{ 2j+1, 2j-1} \mtx S_{2j-1}^{-1} \mtx {\tilde f}_{ 2j-1}$\;
  }
  \For{$j=1,\dots, n_1/b+1$}{
  $\mtx {\tilde f}_{ 2j - 1} \gets \mtx S_{2j-1}^{-1} \mtx {\tilde f}_{ 2j-1}$\;
  }
  \For{$j=n_1/b, \dots, 1$}{
  $\mtx {\tilde f}_{ 2j - 1} \gets \mtx {\tilde f}_{ 2j - 1} - \mtx S_{2j-1}^{-1} \mtx T_{2j-1, 2j+1} \mtx {\tilde f}_{ 2j+1}$\;
  }
  $\mtx {\tilde u} \gets \mtx {\tilde f}$\;
\end{algorithm}

The most elegant way to solve (\ref{eq:Tred_system}) is
to compress the matrices $\mtx S_1, \mtx S_3, \dots$ of Algorithm \ref{alg:sweeping} as HBS using
randomized black-box algorithms; this would lead to a solver with linear complexity in the
case where the PDE is kept fixed as $N$ increases.
However, we have found that for two dimensional problems up to size $N \approx 10^8$, 
the sizes of the slab interfaces are small enough that 
dense linear algebra is the most efficient
way to solve (\ref{eq:Tred_system}) in practice. 
After forming the sub-blocks of $\mtx T$ as HBS, we convert them to dense matrices
as needed in Algorithm \ref{alg:sweeping}. We have observed that
retainining HBS structure in the off-diagonal blocks of $\mtx T$ (as opposed to storing them densely)
is useful to lower the memory footprint of the factorization 
and facilitate fast data movement from the CPU to the GPU in the implementation of  Algorithm \ref{alg:sweeping}.

\vspace{0.5em}

\begin{remark} \label{remark:helm_hbs}
For Helmholtz problems of the form (\ref{eq:var_helm}), we are often interested in evaluating the
performance of a solver where the number of points per wavelength is fixed
as $N$ increases. In this regime, it is still possible to exploit rank structure
in the blocks of $\mtx{T}$; however, the ranks will grow during the execution of 
Algorithm \ref{alg:sweeping}.
To minimize rank-growth, an odd-even elimination order can
be used instead of the sequential one in Algorithm \ref{alg:sweeping}
to keep the slabs as thin as possible during the factorization process.
Due to the substantial rank growth of oscillatory problems in this regime, the HBS algebra
is not likely to perform better than highly-tuned dense linear algebra
routines on the GPU.
Factorizing $\mtx S_1, \mtx S_3, \dots$ of Algorithm \ref{alg:sweeping} using dense linear algebra
allows us to side-step this complication and still achieve excellent performance in Section \ref{sec:hps_numerical},
as our method
relies on the sparsity pattern of the original matrix only and leverages efficient GPU offloading.
\end{remark}

%Though we form $\mtx T$ using efficient black-box algorithms in HBS/HSS matrix format, we factorize $\mtx T$ using dense linear algebra because it is perfectly robust and very fast on a GPU. 
%The HBS/HSS format is useful for reducing the memory requirements of the stored factorization, e.g. the off-diagonal blocks of $\mtx T$ in Algorithm 2 can be stored in HBS/HSS format.
%Algorithm \ref{alg:sweeping} outlines a very simple sweeping algorithm for the factorization, where each block is factorized in serial. 
%Algorithm \ref{alg:sweeping_solve} describes how to apply the factorization to solve systems (\ref{eq:Tred_system}).

\section{Algorithm and complexity costs}

In this section, we provide a summary of the proposed algorithm and discuss its complexity costs.
We also analyze the complexity costs, and choose the buffer size $b$ as a function of the number of grid points $N=n_1 n_2$,
in order to balance costs and achieve competitive complexities for the build and solve times.

We briefly summarize the algorithm for SlabLU.
In Stage One, we compute factorizations of the form
\begin{equation}
\underset{n_2 b \times n_2 b}{\mtx A_{22}^{-1}}, \underset{n_2 b \times n_2 b}{\mtx A_{44}^{-1}}, \dots \label{eq:sparse_blocktri}
\end{equation}
for $n_1/b$ sparse matrices. The reduced matrix $\mtx T$ is constructed using efficient black-box algorithms 
that recover $\mtx T_{11}, \mtx T_{13}, \mtx T_{31}, \mtx T_{33}, \dots$ in HBS format through a random sampling method,
as discussed in Section \ref{sec:hmatrix_rand}. 
It is important to note that Stage One can be trivially parallelized for each slab.

In Stage Two, the reduced system $\mtx T$ is factorized. 
Because the sub-blocks of $\mtx T$ are HBS, this structure can be used create a linear solver. 
However, we have found empirically that for 2D problems, the blocks are small enough that using dense linear algebra is just as fast, and much easier to implement. 
The simplified scheme using dense matrix algebra has an additional benefit in that it relies on sparsity alone; this makes it well suited for highly oscillatory problems, as it is immune to the rank-growth that typically happens when the wavenumber $\kappa$ grows with $N$, cf. Remark \ref{remark:helm_hbs}.

\RestyleAlgo{boxruled}
\begin{algorithm}[!htb]
\caption{Factorizing $\mtx A$ using SlabLU.}
\label{alg:factorA}
\For{j=$1,\dots n_1/b$}{
    \text{Compute $\mtx A_{2j,2j}^{-1}$ using black-box sparse direct solvers.}
}

\text{ \parbox[t]{300pt}{
        Compress $\mtx T_{11}, \mtx T_{13}, \dots$ as HBS matrices
        with randomized black-box compression using matrix-free formulas, e.g. (\ref{eq:T11_detailed}).
        }
}

\text{\parbox[t]{300pt}{
Factorize $\mtx T$ using Algorithm \ref{alg:sweeping}.}
}

\end{algorithm}

\vspace{-1em}

\RestyleAlgo{boxruled}
\begin{algorithm}[!htb]
\caption{Solving $\mtx A \mtx u = \mtx f$ using SlabLU.}
\label{alg:solve}
\text{Calculate the equivalent body $\mtx {\tilde f}$ on $I_1, I_3, \dots$ using (\ref{eq:f1}-\ref{eq:f5}) in parallel.} \hspace{12em}
\tcp{Parallel computation}
\text{Solve $\mtx T \mtx {\tilde u} = \mtx {\tilde f}$ for $\mtx {\tilde u}$ on $I_1, I_3, \dots$ using Algorithm \ref{alg:sweeping_solve}.}
\hspace{12em}\tcp{Serial algorithm}
\text{Solve for $\mtx u$ on $I_2, I_4, \dots$ using $\mtx {\tilde u}$ on $I_1, I_3, \dots$ with (\ref{eq:usol_even}) in parallel. }
\hspace{12em}\tcp{Parallel computation}
\end{algorithm}

\subsection{Algorithmic complexity} 
\label{sec:choosing_b}

When rank-structure is exploited in both Stage One and Stage Two, choosing the slab width $b= \mathcal O(1)$ leads to overall linear complexity (in the case where the equation is kept fixed as $N$ is increased). 
In the simplified scheme relying on dense linear algebra in the second stage, we can choose the buffer size to balance the costs of the two stages.
We first demonstrate how to choose $b$ for the 2nd order finite difference stencil, then repeat the process
for the HPS discretization. The asymptotic costs with $N$ are the same for both discretizations, and as
described in Section \ref{sec:complexity_hps}; there are pre-factors that depend on the local polynomial order $p$ for the HPS discretization.
For both discretizations we report $T_{\rm build}$, which is the asymptotic flop cost to build the direct solver and $M$ which is the memory in bits required to store the direct solver. 
The flop cost to apply the direct solver is $T_{\rm solve}$ and is the same asymptotically as $M$ unless otherwise mentioned.

\subsection{Complexity analysis for 2nd order finite difference discretization}
\label{sec:complexity_fd}

Stage One requires computing sparse factorizations (\ref{eq:sparse_blocktri}) for each thin slab,
then constructing the subblocks  of $\mtx T$ in HBS format using a randomized black-box algorithm.
For each slab, the local sparse matrix is of size $n_2 b \times n_2 b$ with bandwidth $b$; the sparse factorization
requires $\mathcal O (n_2 b^3)$ flops to compute and $\mathcal O (n_2 b^2)$ bits to store.
Each sub-block of the reduced system $\mtx T$ is constructed as an HBS matrix
using a randomized sampling algorithm, for which the dominant cost is constructing the random samples
 (\ref{eq:rand_samples}), which requires $\mathcal O\left(n_2 b^3 \right)$ flops using matrix-free formulas for applying (cf. equation (\ref{eq:T11_detailed}) for the formula to apply $\mtx T_{11}$).

Stage Two involves a factorization of $\mtx T$ using Algorithm \ref{alg:sweeping}.
For generality, consider that $T_{ \mtx S}$ and $M_{ \mtx S}$ are the time and memory, respectively
needed for the factorization of each matrix $\mtx S_1, \mtx S_3, \dots$ of size $n_2 \times n_2$. Then the total costs
of computing and storing SlabLU are
\begin{equation}
T_{\rm build} = \underset{\substack{\text{sparse slab}\\ \text{factorizations}}}{\mathcal O \left( \frac{n_1} b n_2 b^3 \right)}
 + \underset{\text{factorization of $\mtx T$}}{ \mathcal O \left (\frac {n_1} b T_{\mtx S} \right) },
\qquad M = \underset{\substack{\text{sparse slab}\\ \text{factorizations}}}{\mathcal O \left( \frac{n_1} b n_2 b^2 \right)}
 + \underset{\text{factorization of $\mtx T$}}{ \mathcal O \left (\frac {n_1} b M_{\mtx S} \right) }.
\label{eq:complexity_funcb}
\end{equation}
Because direct solvers are typically 
limited by their memory costs, we choose $b$ to miminize the memory footprint. This depends on the choice of matrix format used for the matrices $\mtx S_1, \mtx S_3, \dots$ 
For instance, the use of HBS algebra gives $T_{\mtx S} = \mathcal O \left (k\ n_2 \right)$ and
 $M_{\mtx S} =\mathcal O \left( k^2\ n_2 \right)$. When the PDE is fixed as $N$ grows, the resulting complexity is
\begin{equation}
b = \mathcal O \left(\sqrt {k} \right)\qquad  \Rightarrow \qquad T_{\rm build} 
= \mathcal O \left(  k^{1.5} N \right), \qquad M = \mathcal O \left(  k^{0.5} N \right)
\label{eq:complexity_hbs}
\end{equation}
The use of dense linear algebra to factorize $\mtx T$ gives
\begin{equation}
b = \mathcal O \left(\sqrt {n_2} \right)\qquad \Rightarrow \qquad 
T_{\rm build} = {\mathcal O \left( n_1 n_2^{2.5} \right)},
\qquad M = {\mathcal O \left( n_1 n_2^{1.5} \right)},
\label{eq:complexity_dense}
\end{equation}
which scales particularly well for domains with high aspect ratio when $n_1 > n_2$. 
SlabLU provides a flexible framework for choosing $b$ depending on the choice of rank-structured linear algebra in factorizing $\mtx T$. Likewise, the buffer size $b$ can be chosen for block low-rank formats \cite{amestoy2017complexity}.
The convenience of SlabLU lies in the fact that the matrices $\mtx S_1,\mtx S_3, \dots$ are roughly the same size, and the $\mathcal H$ or $\mathcal H^2$ algebra is much simpler to optimize for performance.

In our numerical results, we used dense linear algebra to factorize $\mtx T$, which is highly
effective and robust for oscillatory problems in 2D, up to $N \approx 10^8$. Section \ref{sec:num} features a variety of rectangular and curved domains, as well as square domains, which is the adversarial worst case in terms of the algorithm complexity in (\ref{eq:complexity_dense}). When $n_1=n_2$, the complexity is
$T_{\rm build} = \mathcal O\left( N^{1.75} \right)$ and $ M = \mathcal O\left( N^{1.25} \right)$,
though complexities observed in Section \ref{sec:num} for $T_{\rm build}$ are practically linear, especially 
when using GPU off-loading.

\subsection{Complexity analysis for HPS discretization}
\label{sec:complexity_hps}

The numerical results feature a high order discretization scheme that interfaces particularly well with sparse direct solvers
which we use to resolve high frequency scattering problems to high accuracy.
As discussed in Section \ref{sec:other_disc}, HPS is a spectral collocation discretization that employs multiple subdomains
to enforce the PDE using spectral differentiation. The subdomains are coupled together by ensuring continuity of the solution 
and its derivative across subdomains. 

A natural approach to factorizing the sparse coefficient matrix is to first factorize each leaf subdomain in parallel
with $\mathcal O\left({p^4 N}\right)$ flops, cf. Figure \ref{fig:hps_two_sub} for a description of static condensation. 
Then, the remaining sparse matrix $\mtx {\tilde A}$ of size roughly $N/p \times N/p$ is factorized using SlabLU.
Because $\mtx{\tilde A}$ is smaller, the costs of factorizing $\mtx {\tilde A}$ with SlabLU involves prefactors with $p$.
In particular, the sparse factorization of a thin slab of width $b$ requires $\mathcal O \left({\frac{n_2}{\sqrt{p}} b^3} \right)$ flops to compute and $\mathcal O \left({\frac{n_2}{\sqrt{p}} b^2} \right)$ bits to store.
As we did in equation (\ref{eq:complexity_funcb}), we can choose $b$ to minimize the memory cost of storing the factorization of $\mtx {\tilde A}$, leading to prefactors of $\sim 1/\sqrt{p}$ in equations (\ref{eq:complexity_hbs}, \ref{eq:complexity_dense}).

Traditionally, the pre-factor cost of static condensation has been considered prohibitively expensive. 
However, \cite{yesypenko2022parallel} describes simple GPU optimizations that use batched linear algebra to significantly accelerate these operations; they are so efficient that we can save on storage costs by not explicitly storing
the factorizations of the local spectral differentiation matrices in each HPS subdomain. Instead, 
we can reform and refactor these matrices as needed, leading to a decreased cost in storage (only the factorization of $\mtx {\tilde A}$ is stored)
and an increased cost in applying the factorization of $\mtx A$. 

\section{Numerical experiments}
\label{sec:num}

In this section, we demonstrate the effectiveness of our solver through a series of numerical experiments.
We report the build time, solve time, and accuracy of solving constant and variable-coefficient elliptic PDEs using
two collocation-based discretization schemes on both rectangular and curved geometries. Our experiments were conducted
on various hardware architectures to showcase the portability and ease of performance tuning of our framework.
We utilize a high-order multidomain spectral collocation scheme, briefly introduced in Section \ref{sec:other_disc},
to solve challenging scattering phenomena. The high-order discretization scheme allows us to accurately discretize the PDE, 
while the flexibility of the multidomain scheme enables us to solve PDEs on curved domains using SlabLU. The combination of SlabLU
and high-order discretization provides a powerful tool for simulating electromagnetic and acoustic scattering. 

We have implemented SlabLU in an open-source software package, \texttt{slabLU}, in Python \cite{YesypenkoSlabLUATwoLevel2024}. It is designed to be efficient and portable across various hardware architectures. The source code, detailed documentation, and usage examples are provided to facilitate the reproduction of our numerical results.
We conducted the experiments on three different architectures: (1) a 16-core Intel i9-12900k CPU with 128 GB of RAM,
(2) an NVIDIA RTX 3090 GPU with 24 GB of memory and access to 128 GB of RAM, and (3) an NVIDIA V100 with 32 GB of memory
and access to 768 GB of RAM. We chose to run experiments on architectures (1) and (2) to demonstrate that the memory volume required
to run SlabLU is reasonable. All experiments used double precision.

In our implementation, we use GPU offloading for Stage One and Stage Two in the factorization of $\mtx A$ using SlabLU, for architectures (2) and (3).
In Stage One, we wrote custom sparse direct solvers for sparse systems (\ref{eq:sparse_blocktri}) using GPU acceleration.
This allowed us to achieve substantial acceleration in generating random sketches (\ref{eq:rand_samples}) of each subblock of $\mtx T$.
The implementation also uses batched linear algebra \cite{abdelfattah2021set} to post-process sketches and recover the subblocks of $\mtx T$
as HBS matrices. In Stage Two, we use the HBS structure of $\mtx T$ to efficiently move the reduced
system onto the GPU, then the matrices $\mtx S_1, \mtx S_3, \dots$ of Algorithm  \ref{alg:sweeping} are factorized 
densely using highly optimized vendor libraries on the GPU.
Recent works \cite{ghysels2022high,abdelfattah2022addressing} have demonstrated promising results in achieving high performance for sparse direct solvers 
on the GPU. Established library packages, such as STRUMPACK \cite{strumpack_doecode}, have also incorporated GPU support. 
The challenge in achieving optimal performance for such codes lies in careful load balancing and memory access patterns across a range of front sizes
in a multi-level tree \cite{vuduc2010limits,kim2013scheduling}.
SlabLU circumvents this difficulty using a two-level scheme that uses traditional sparse direct solvers for small fronts of size $b$ and
custom $\mathcal H^2$-matrix approaches for the larger fronts.

\subsection{Description of benchmark PDEs and accuracies reported}
We describe the PDEs used as benchmarks in our numerical experiments and how we calculate accuracy.
We use the constant coefficient Helmholtz problem for various wavenumbers $\kappa$
\begin{equation}
\left\{\begin{aligned}
-\Delta u(x)  - \kappa^2 u(x) =&\ 0,\qquad&x \in \Omega,\\
u(x) =&\ u_{\text{true}}(x),\qquad&x \in \Gamma,
\end{aligned}\right. 
\label{eq:constant_coeff_helm}
\end{equation}
where the true solution $u_{\text{true}} = J_0 \left( \kappa \|x - (-0.1,0.5)\| \right)$
and $J_{0}$ is the zeroth Bessel function of the first kind.

After applying the direct solver, we obtain the calculated solution $\mtx u_{\rm calc}$ at discretization points within the domain.
We report the relative error with respect to the residual of the discretized system (\ref{eq:Au=f}) and with respect to the
true solution $\mtx u_{\rm true}$ evaluated at the collocation points as follows:
\begin{equation} {\rm relerr}_{\rm res} = \frac{{\|\mtx A \mtx u_{\text{calc}} - \mtx f\|}_2}{{\|\mtx f\|}_2}, 
\qquad {\rm relerr}_{ \rm true} = \frac{{\|\mtx u_{\rm calc} - \mtx u_{\rm true}\|}_2}{{\| \mtx u_{\rm true}\|}_2}.\end{equation}

We also report $T_{\rm build}$, which is the wall-clock time needed to factorize the coefficient matrix $\mtx A$, and $M$, 
which is the memory required to store the computed direct solver. Additionally, we report $T_{\rm solve}$ for one right-hand side vector, 
which is the time needed to apply the direct solver to solve systems (\ref{eq:Au=f}), as described in Algorithm \ref{alg:solve}.

\subsection{Experiments using Low-Order Discretization}

In this section, we demonstrate how SlabLU performs on sparse coefficient matrices arising from PDEs discretized with 2nd order finite differences.
We also compare SlabLU to SuperLU, a black-box sparse direct solver code.

\begin{figure}[!htb]
\begin{subfigure}{\textwidth}
\centering
\includegraphics[width=0.65\textwidth]{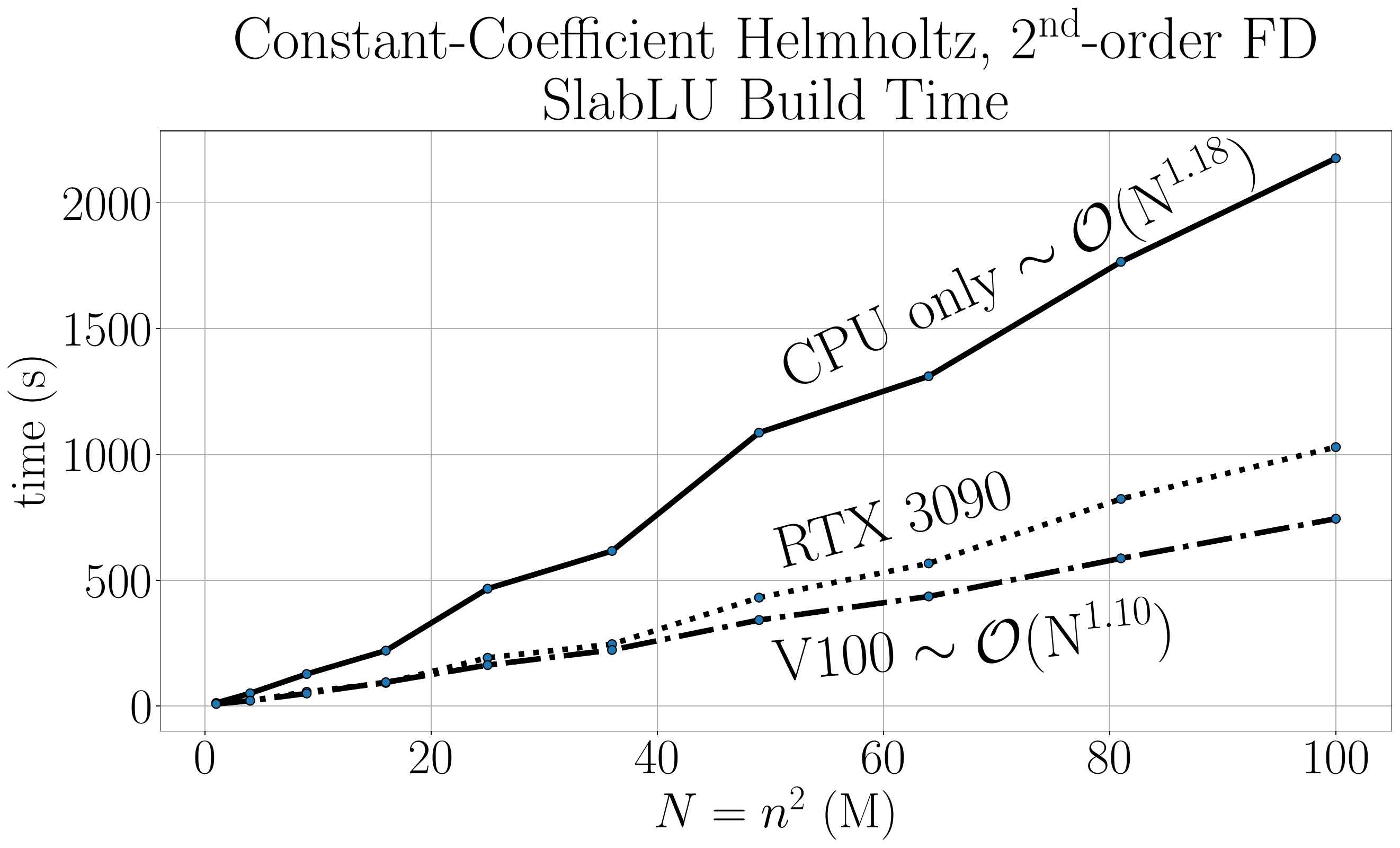}
\caption{Timing results for factorizing $\mtx A$ using SlabLU on various architectures.}
\end{subfigure}
\begin{subfigure}{\textwidth}
\centering
\scriptsize
\begin{tabular}{rr|r|rr|rr}
$N$ & $b$ & $\kappa$ & $M$ &$T_{\rm solve}$ & ${\rm relerr}_{\rm res}$ & ${\rm relerr}_{\rm true}$\\ \hline
1.0 M & 50 & 27.1 & 0.5 GB & 0.9 s & 1.1e-11 & 1.8e-03 \\
4.0 M & 100 & 52.3 & 2.5 GB & 2.9 s & 1.8e-11 & 2.6e-03 \\
9.0 M & 125 & 77.4 & 6.6 GB & 7.3 s & 2.6e-11 & 4.0e-03 \\
16.0 M & 160 & 102.5 & 12.9 GB & 10.4 s & 2.4e-11 & 1.0e-02 \\
25.0 M & 200 & 127.7 & 22.5 GB & 15.2 s & 2.5e-11 & 1.0e-02 \\
36.0 M & 200 & 152.8 & 32.3 GB & 21.7 s & 2.8e-11 & 4.8e-03 \\
49.0 M & 200 & 177.9 & 46.9 GB & 31.8 s & 4.2e-11 & 6.5e-03 \\
64.0 M & 250 & 203.1 & 64.0 GB & 42.8 s & 1.4e-10 & 3.1e-02 \\
81.0 M & 250 & 228.2 & 83.8 GB & 75.5 s & 8.5e-11 & 1.1e-02 \\
100.0 M & 250 & 253.3 & 105.7 GB & 88.6 s & 6.7e-11 & 8.7e-03 \\
        &     & & $\sim \mathcal O(N^{1.16})$ & $\sim \mathcal O(N^{1.09})$&\\
\end{tabular}
\caption{\small The table reports how $b$ grows as a function of the problem size $N$, as well as $M, T_{\rm solve}$
and the computed relative accuracies. The wavenumber $\kappa$ is increased with the problem size to maintain 250 points per wavelength. Though the solution is resolved to at least 10 digits in the residual, the relative error compared to the true solution is only 3 digits. For these experiments, $T_{\rm solve}$ is reported on the CPU architecture.}
\end{subfigure}
\caption{Helmholtz equation (eq. \ref{eq:constant_coeff_helm}) discretized with 2nd order finite differences with constant points per wavelength.}
\label{fig:greens_fd}
\end{figure}

We demonstrate competitive scaling for the build time of the factorization and for the memory footprint. 
See Figure \ref{fig:greens_fd} for the constant-coefficient Helmholtz equation (\ref{eq:constant_coeff_helm}).
Despite the super-linear complexity scaling, the scaling appears to be linear for grids of size up to $N=100$M.
We discretize the Helmholtz equation to at least 10 points per wavelength to resolve the oscillatory solutions,
leading to large grids. However, due to the effect of pollution when using low-order discretization, 
we need to discretize the Helmholtz equation to 250 points per wavelength
to attain 3 digits of accuracy with respect to the known analytic solution.

We compare the performance of SlabLU and SuperLU in solving the 2nd order finite difference discretization of the constant coefficient Helmholtz
equation (\ref{eq:constant_coeff_helm}), which results in very ill-conditioned sparse matrices (\ref{eq:Au=f}) that need to be solved.
SuperLU \cite{bollhofer2020state,2011_li_supernodal} is a generic sparse direct solver that computes a sparse LU factorization of any 
given sparse matrix to high accuracy.
We use the Scipy interface (version 1.8.1) to call SuperLU with the default permutation specification of COLAMD ordering; this
version of SuperLU is sequential with multi-threading for BLAS calls.

\begin{figure}[!htb]
\begin{subfigure}{0.55\textwidth}
\centering
\includegraphics[width=\textwidth]{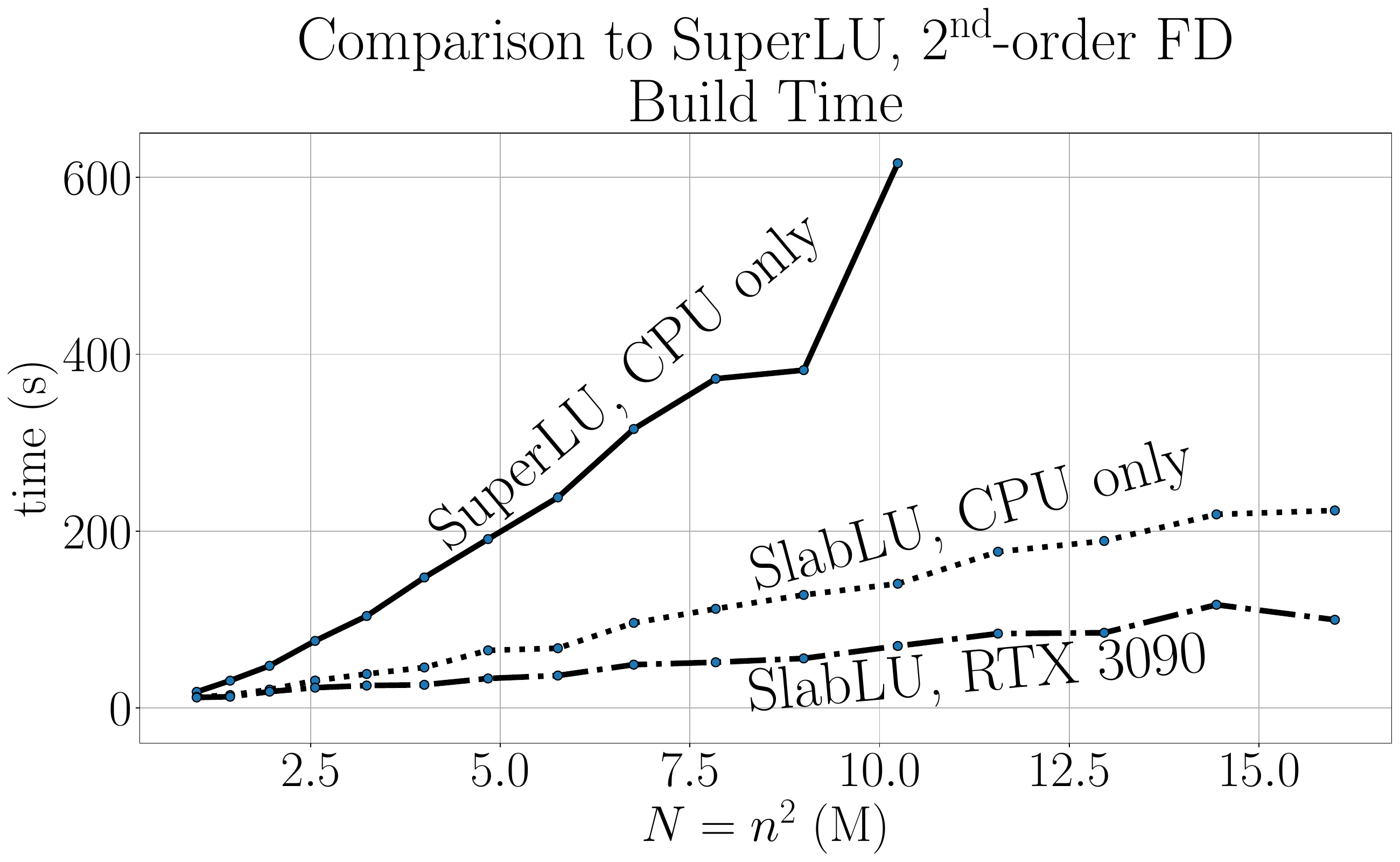}
\caption{Time to factorize $\mtx A$.}
\end{subfigure}%
\hfill
\begin{subfigure}{0.45\textwidth}
\scriptsize
\centering
\begin{tabular}{rr|rr}
 & & \multicolumn{2}{c}{$T_{\rm solve}$}\\ \hline
$N$ & $\kappa$ & SlabLU & SuperLU \\ \hline
1.0 M & 27.12 & 0.21 s& 0.46 s\\
1.4 M & 27.12 & 0.31 s& 0.51 s\\
2.0 M & 33.40 & 0.43 s& 0.69 s\\
2.6 M & 39.69 & 0.58 s& 1.03 s\\
3.2 M & 45.97 & 0.74 s& 1.27 s\\
4.0 M & 52.25 & 0.69 s& 1.38 s\\
4.8 M & 52.25 & 1.16 s& 1.92 s\\
5.8 M & 58.54 & 1.10 s& 2.25 s\\
6.8 M & 64.82 & 1.69 s& 2.48 s\\
7.8 M & 71.10 & 1.96 s& 3.21 s\\
9.0 M & 77.39 & 1.63 s& 3.47 s\\
10.2 M & 77.39 & 1.97 s& 3.84 s\\
11.6 M & 83.67 & --- & 4.78 s
\end{tabular}
\caption{Table reporting $\kappa$ and $T_{\rm solve}$ on the CPU architecture.}
\end{subfigure}

\vspace{1em}
\begin{subfigure}{0.49\textwidth}
\centering
\includegraphics[width=0.85\textwidth]{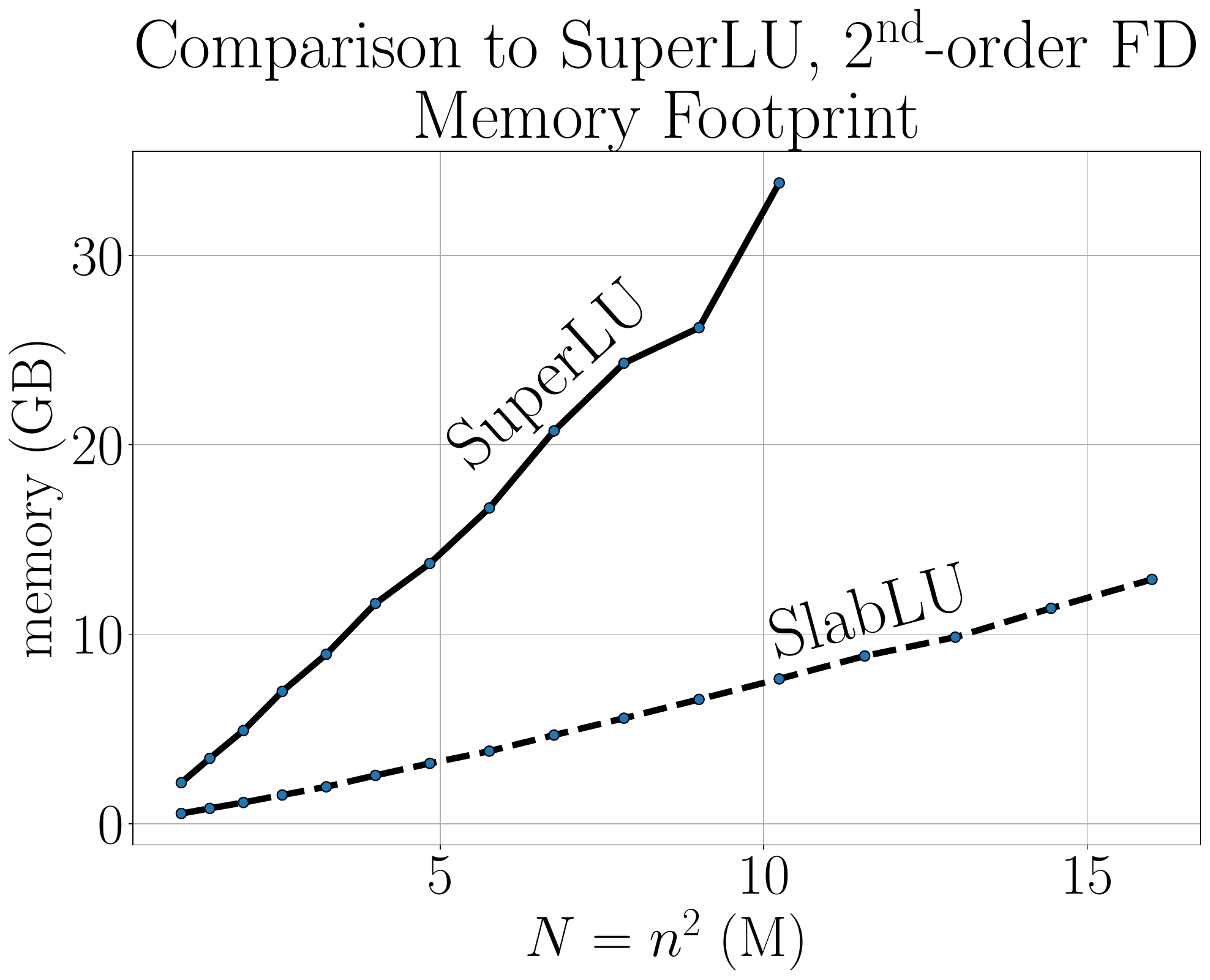}
\caption{Memory need to store factorization of $\mtx A$.}
\end{subfigure}%
\hfill
\begin{subfigure}{0.49\textwidth}
\centering
\includegraphics[width=0.9\textwidth]{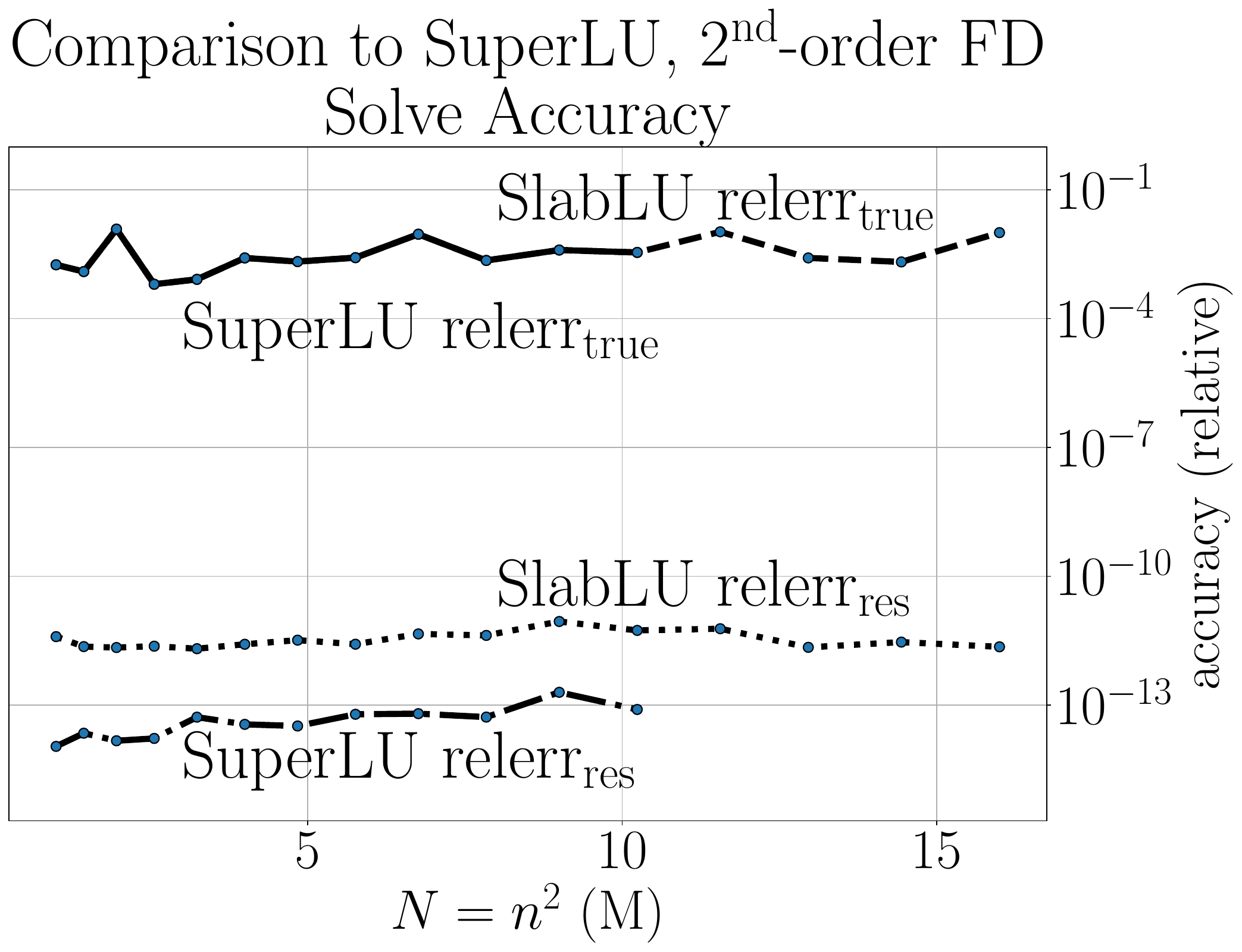}
\caption{Relative accuracies of computed solutions.}
\end{subfigure}
\caption{\small Comparison to a sequential version of Super LU for 2d order FD which uses multi-threading through BLAS calls.
SlabLU uses approximations (cf. Algorithm \ref{alg:factorA}) for construction of the reduced system $\mtx T$
and efficiently factorizes the system (\ref{eq:Au=f}) to high accuracy.
Our comparison shows that SlabLU outperforms this version of SuperLU in terms of build times and memory costs.
For the CPU-comparison, SlabLU is faster by a factor of 4 for $N$=10.2M. Using GPU acceleration
makes the method faster by a factor of 8. Additionally, SlabLU is more memory efficient by a factor of 4 for $N = 10.2 \rm{M}$.}
\label{fig:SuperLU_fd}
\end{figure}

The primary purpose of the comparison to this version of SuperLU is to compare accuracy of the computed factorization.
Sparse direct solvers use sparsity in the discretized operator in order to factorize the
sparse coefficient matrix exactly. SlabLU uses sparsity in the traditional sense for the sparse
factorizations (\ref{eq:sparse_blocktri}) in Stage One and approximately constructs the reduced system $\mtx T$
to high accuracy using randomized black-box sampling (cf. Section \ref{sec:hmatrix_rand}.)
Both schemes are able to resolve the solution up to the discretization error (cf. Figure \ref{fig:SuperLU_fd}).

SlabLU also compares favorably to SuperLU in terms of memory costs; this is slightly surprising because the
memory costs of SlabLU have slighly worse asymptotic scaling, compared to multi-level schemes. 
We believe may be because of inefficiencies of storing the sparse factorization or because COLAMD chooses a suboptimal ordering, compared to a METIS ordering, which is not available for this version. 
For problems larger than 10.2M points, SuperLU does not compute the factorization, because the memory requirements exceed some pre-prescribed limit.
The authors note that more optimized solvers are available that use multi-threading or MPI parallelism \cite{li2003superlu_dist,2004_davis_UMFPACK,amestoy2000mumps}; 
comparison to these solvers will be the subject of future work.

\subsection{Solving challenging scattering problems with high order discretization}
\label{sec:hps_numerical}

\begin{figure}[!htb]
\begin{subfigure}{\textwidth}
\centering
\includegraphics[width=0.65\textwidth]{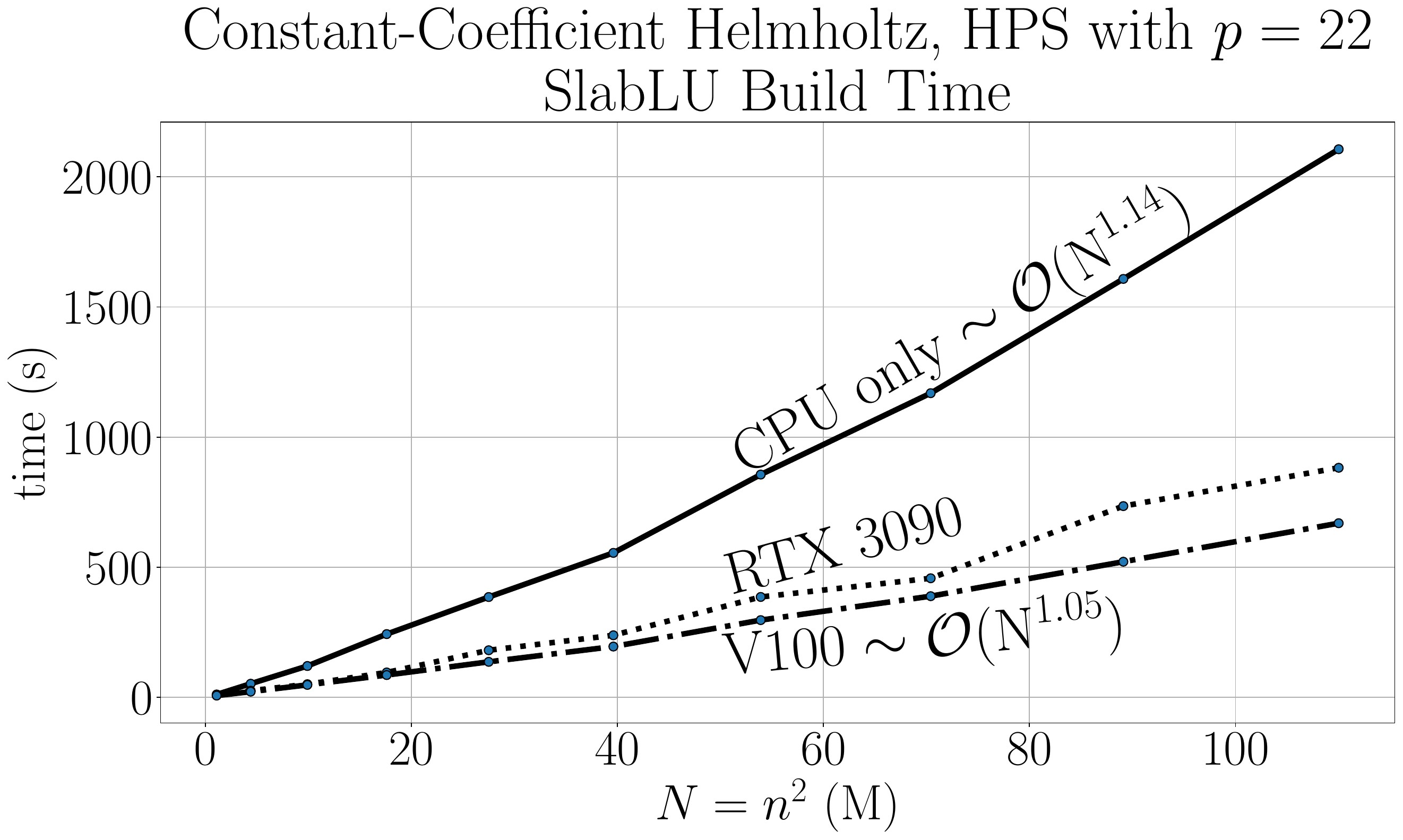}
\caption{Timing results for factorizing $\mtx A$ using SlabLU on various architectures.}
\end{subfigure}
\begin{subfigure}{\textwidth}
\centering
\scriptsize
\begin{tabular}{rr|r|rr|rr}
$N$ & $b$ & $\kappa$ & $M_{\rm build}$& $T_{\rm solve}$ & $\text{relerr}_{\text{res}}$ & $\text{relerr}_{\text{true}}$\\ \hline
1.1 M & 100 & 630.3 & 0.3 GB & 1.0 s & 4.2e-12 & 2.4e-08 \\
4.4 M & 200 & 1258.6 & 1.8 GB & 3.8 s & 6.5e-12 & 9.5e-08 \\
9.9 M & 200 & 1886.9 & 4.4 GB & 7.7 s & 1.2e-11 & 2.2e-07 \\
17.6 M & 200 & 2515.3 & 8.9 GB & 13.9 s & 2.7e-11 & 4.0e-07 \\
27.5 M & 200 & 3143.6 & 15.6 GB & 19.8 s & 6.4e-11 & 9.4e-07 \\
39.6 M & 300 & 3771.9 & 21.5 GB & 33.2 s & 4.3e-10 & 7.3e-07 \\
53.9 M & 280 & 4400.2 & 32.3 GB & 46.4 s & 6.9e-11 & 7.9e-07 \\
70.4 M & 400 & 5028.5 & 45.8 GB & 60.9 s & 1.3e-09 & 5.0e-07 \\
89.1 M & 360 & 5656.9 & 61.0 GB & 76.6 s & 4.3e-10 & 1.2e-06 \\
110.0 M & 400 & 6285.2 & 80.0 GB & 83.9 s & 2.3e-10 & 7.5e-07 \\
        &     & & $\sim \mathcal O(N^{1.17})$ & $\sim \mathcal O(N^{1.00})$&\\
\end{tabular}
\caption{The table reports how $b$ grows as a function of the problem size $N$, as well as $M_{\rm factor}, T_{\rm solve}$
and the computed relative accuracies. The wavenumber $\kappa$ is increased with the problem size to maintain 10 points per wavelength. Using a high order multidomain spectral collocation scheme allows us to avoid the effects of pollution and achieve at least 6 digits of relative accuracy, compared to the known solution.
For these experiments, $T_{\rm solve}$ is reported on the NVIDIA V100 architecture.}
\label{tab:greens_hps}
\end{subfigure}
\caption{Helmholtz equation (eq. \ref{eq:constant_coeff_helm}) discretized with HPS discretization ($p=22$) with constant points per wavelength.}
\label{fig:greens_hps}
\end{figure}

High-order discretization is crucial for resolving variable-coefficient scattering phenomena due to the pollution effect, 
which requires increasing the number of points per wavelength as the domain size increases \cite{beriot2016efficient,deraemaeker1999dispersion}. 
The HPS discretization (cf.~Section \ref{sec:other_disc}) is less sensitive to pollution because
it allows for a high choice of local polynomial order \cite{gillman2014direct,martinsson2013direct}. In Helmholtz experiments, HPS with $p=22$ accurately resolves oscillatory solutions 
using only 10 points per wavelength on domains up to size $1000\lambda \times 1000\lambda$, while 2nd order FD requires 100-250 points per wavelength to achieve low accuracy. Combining HPS with SlabLU provides a powerful tool for resolving challenging scattering phenomena to high accuracy, especially for situations where 
efficient preconditioners are not available \cite{ernst2012difficult}.
We first demonstrate the performance of SlabLU for sparse coefficient matrices arising from the HPS discretization
on a benchmark PDE of constant coefficient Helmholtz in Figure \ref{fig:greens_hps}. As discussed in Section \ref{sec:complexity_hps}, the leaf operations are handled efficiently using batched linear algebra, and the local leaf factorizations are discarded and re-factorized as needed during the solve stage to save on memory costs for the direct solver.

We will now demonstrate the ability of HPS, combined with SlabLU as a sparse direct solver, to solve complex scattering phenomena
on various 2D domains. For the presented PDEs, we will show how the accuracy of the calculated solution converges to a reference 
solution depending on the choice of $p$ in the discretization. Specifically, we will solve the BVP (\ref{eq:bvp}) 
with the variable-coefficient Helmholtz operator (\ref{eq:var_helm}) for Dirichlet data on curved and rectangular domains.

\begin{figure}[!htb]
    \begin{subfigure}{\textwidth}
    \includegraphics[width=\textwidth]{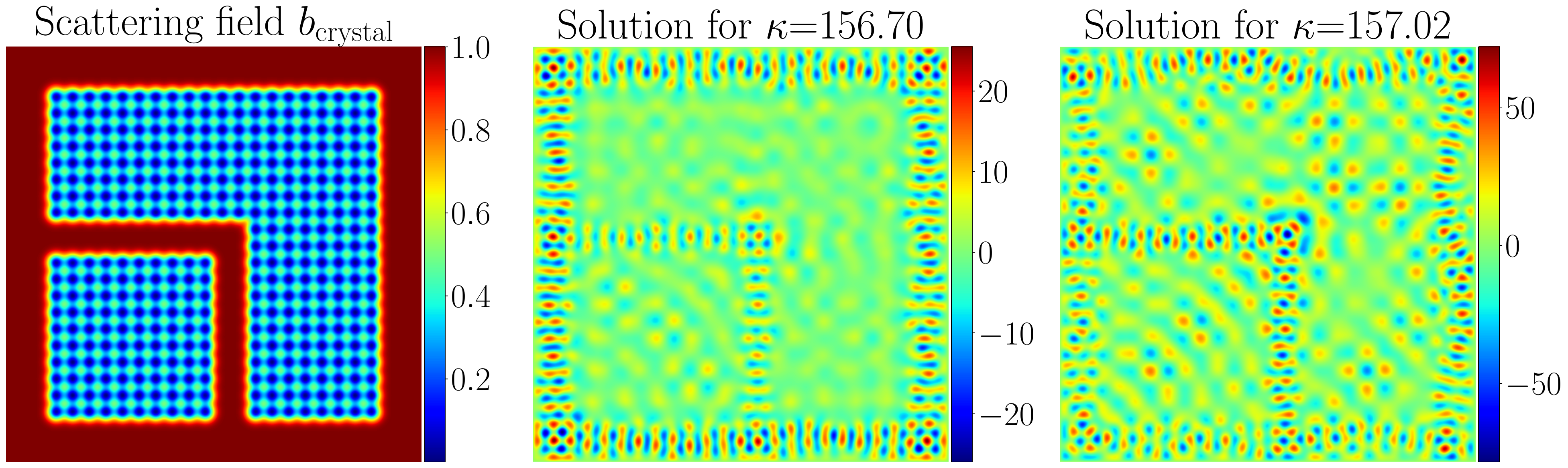}
    \caption{Solutions on $\Omega$.}
    \vspace{-1em}
    \end{subfigure}
    \begin{subfigure}{\textwidth}
    \begin{minipage}{0.5\textwidth}
    \centering
    \vspace{1em}
    \includegraphics[width=0.85\textwidth]{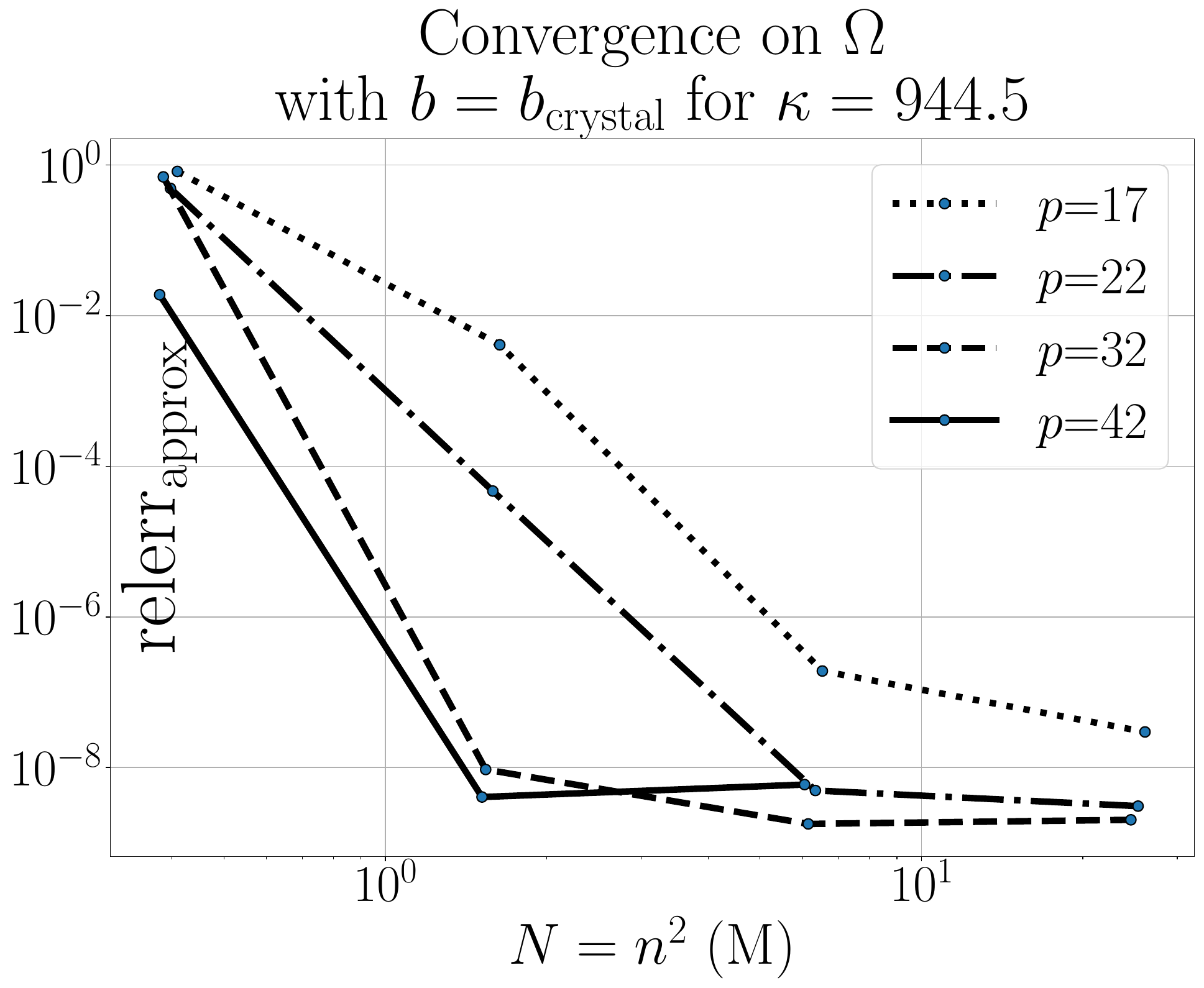}
    \end{minipage}%
    \hfill
    \begin{minipage}{0.45\textwidth}
    \caption{Convergence on square domain $\Omega$ for reference solution $\mtx u_{\rm ref}$ on HPS discretization for $N$=36M with $p=42$.}
    \label{fig:convergence_square}
    \end{minipage}
    \end{subfigure}
    \caption{ Solutions of variable-coefficient Helmholtz problem on domain $\Omega$ with Dirichlet data $u \equiv 1$ on $\partial \Omega$ for various wavenumbers $\kappa$. The scattering field $b_{\rm crystal}$ is a photonic crystal with an extruded corner waveguide. The crystal is represented as a series of narrow Gaussian bumps with separation $s=0.04$ and is designed to filter wave frequencies that are roughly $1/s$.}
    \label{fig:crystal_rhombus}
\end{figure}

We fix the PDE and refine the mesh to compare calculated solutions to a reference solution obtained on a fine mesh with high $p$,
as the exact solution is unknown. The relative error is calculated by comparing $\mtx u_{\rm calc}$ to the reference solution $\mtx u_{\rm ref}$ at a small number of collocation points $\{ x_j \}_{j=1}^M$ using the $l_2$ norm
\begin{equation}\rm {relerr}_{\rm approx} = \frac{{\| \mtx u_{\rm calc} - \mtx u_{\rm ref}\|}_2}{{\|\mtx u_{\rm ref}\|}_2}.\end{equation}
We demonstrate the convergence on a unit square domain $\Omega = [0,1]^2$ with a variable coefficient field $b_{\rm crystal}$ corresponding to a photonic crystal, shown in Figure \ref{fig:crystal_rhombus}.

\begin{figure}[!htb]
\begin{subfigure}{0.55\textwidth}
\centering
\includegraphics[width=0.95\textwidth]{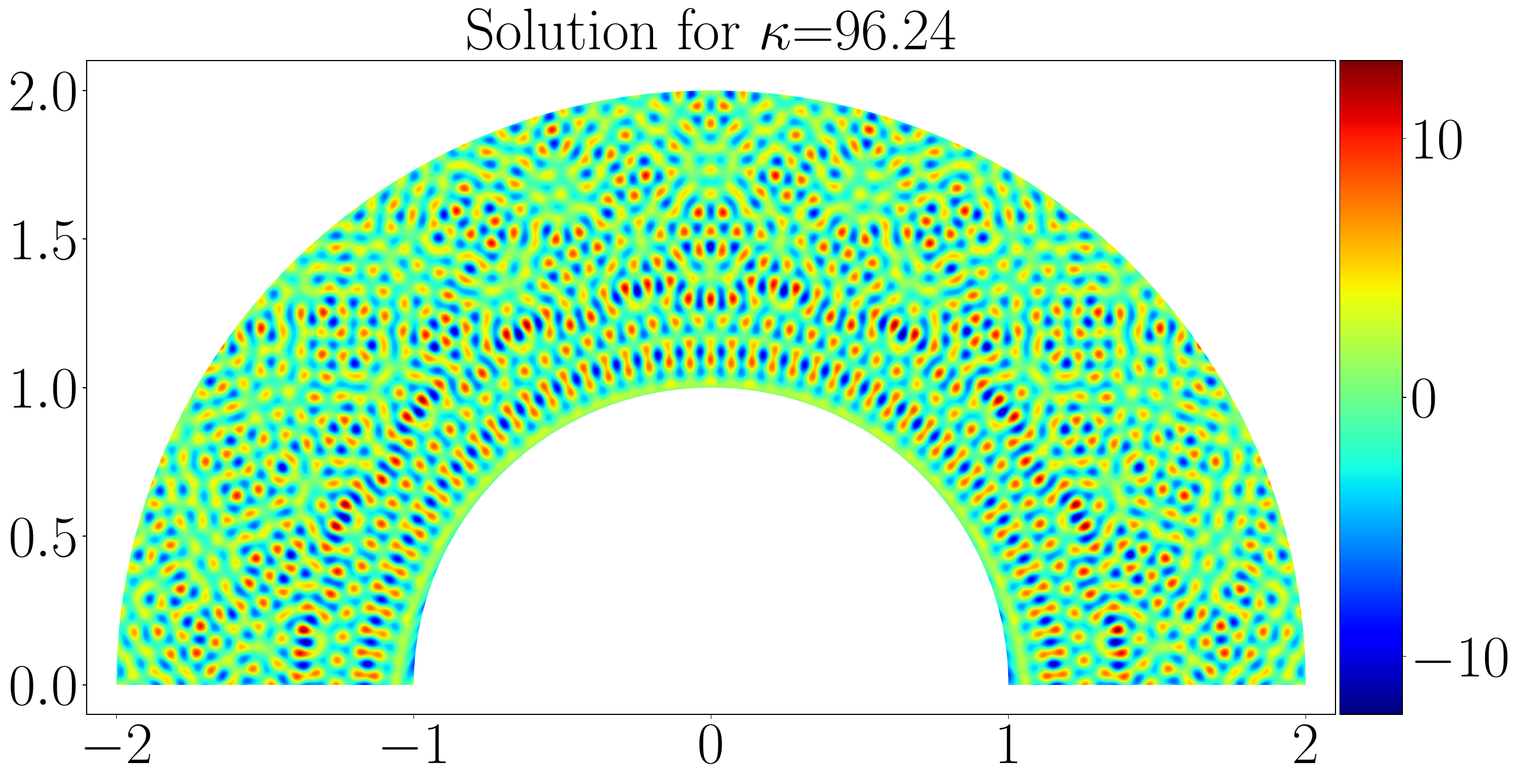}
\vspace{3em}
\caption{Solution on $\Psi$.}
\end{subfigure}%
\hfill
\begin{subfigure}{0.45\textwidth}
\centering
\includegraphics[width=0.95\textwidth]{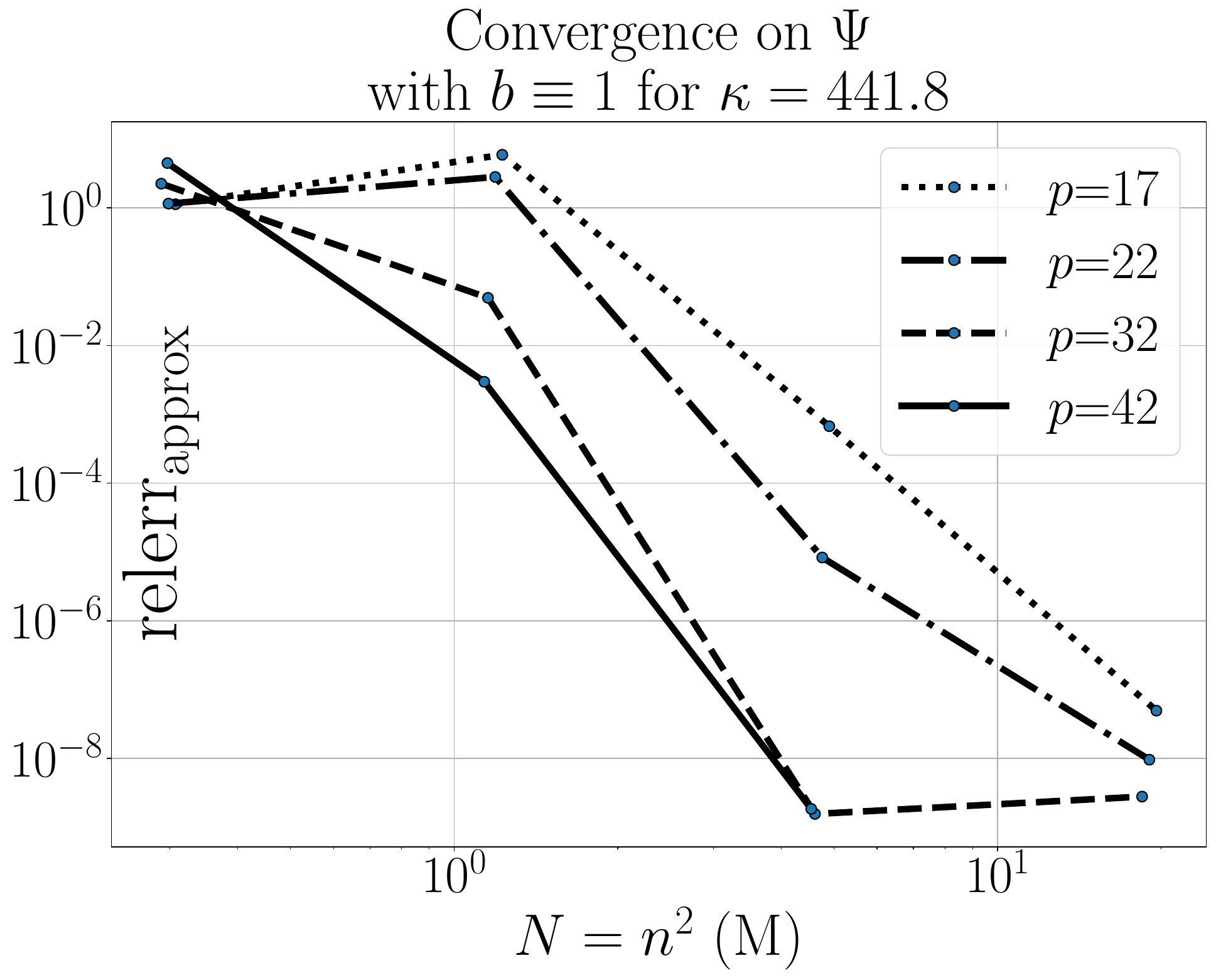}
\caption{Convergence on curved domain $\Psi$ for reference solution $\mtx u_{\rm ref}$ on HPS discretization for $N$=36M with $p=42$.}
\end{subfigure}
\caption{Solutions of constant-coefficient Helmholtz problem on curved domain $\Psi$ with Dirichlet data given by $u \equiv 1$ on $\partial \Psi$. The solutions are calculated parametrizing $\Psi$ in terms of a reference rectangle domain as (\ref{eq:parametrization_annulus}) and solving a variable-coefficient elliptic PDE on the reference rectangle.}
\label{fig:sol_annulus}
\end{figure}

Next, we show the convergence on a curved domain $\Psi$ with a constant-coefficient field 
$b \equiv 1$, where $\Psi$ is a half-annulus given by an analytic parameterization over a reference rectangle. 
The domain $\Psi$ is parametrized as
\begin{equation}
    \Psi =\left\{ \left( \cos \left( \hat{\theta}(x_1) \right) , \sin \left( \hat{\theta}(x_1) \right) \right)\ \text{for}\ (x_1,x_2) \in [0,3] \times [0,1] \right\},
\label{eq:parametrization_annulus}
\end{equation}
where $\hat{\theta}(z) = \frac \pi 3 z$.
Using the chain rule, (\ref{eq:var_helm}) on $\Psi$ takes a different form of a variable-coefficient elliptic PDE on the reference rectangle.
The solutions on $\Psi$ and convergence plot are shown in Figure \ref{fig:sol_annulus}.

\begin{figure}[!htb]
\begin{subfigure}{\textwidth}
\centering
\includegraphics[width=.95\textwidth]{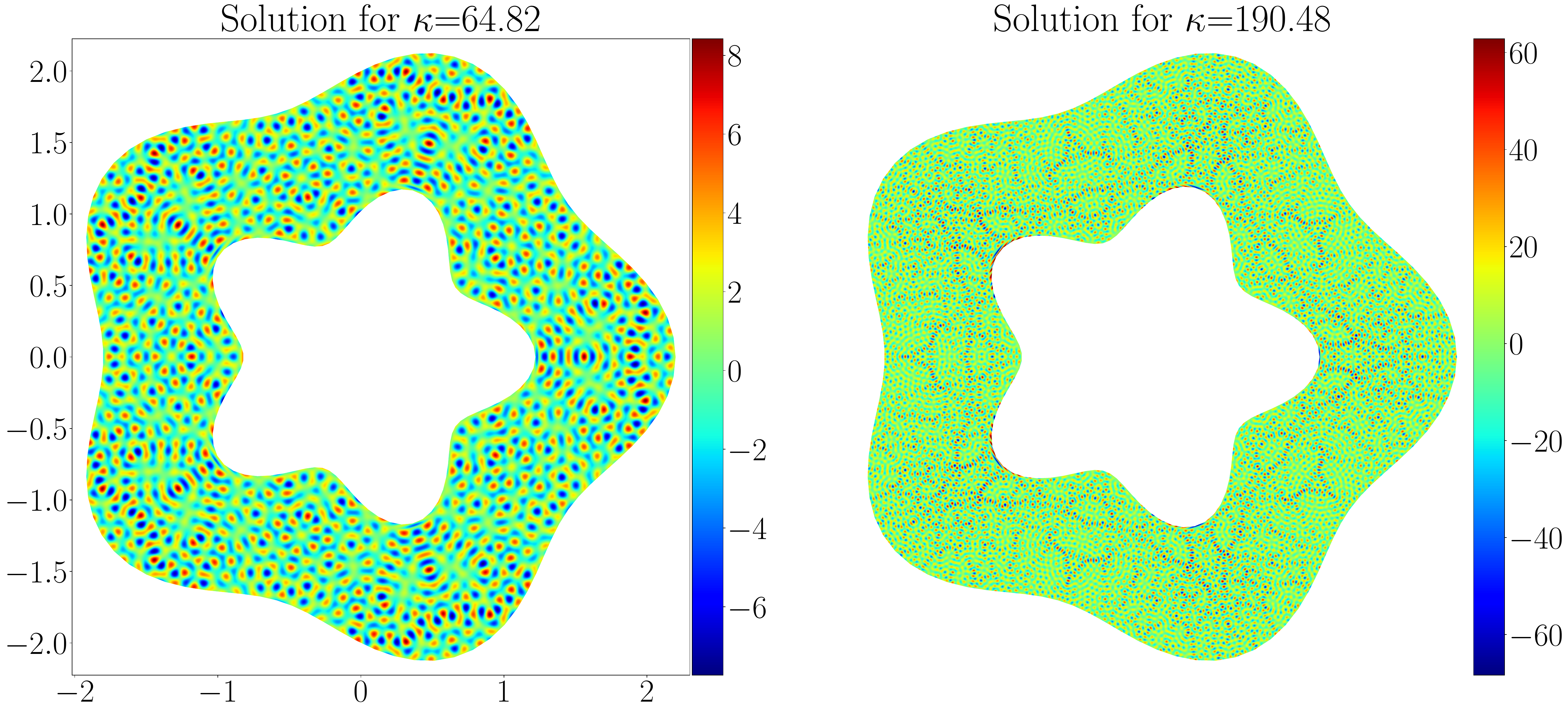}
\caption{Solutions on $\Phi$.}
\end{subfigure}

\begin{subfigure}{\textwidth}
\begin{minipage}{0.5\textwidth}
\centering
\includegraphics[width=0.85\textwidth]{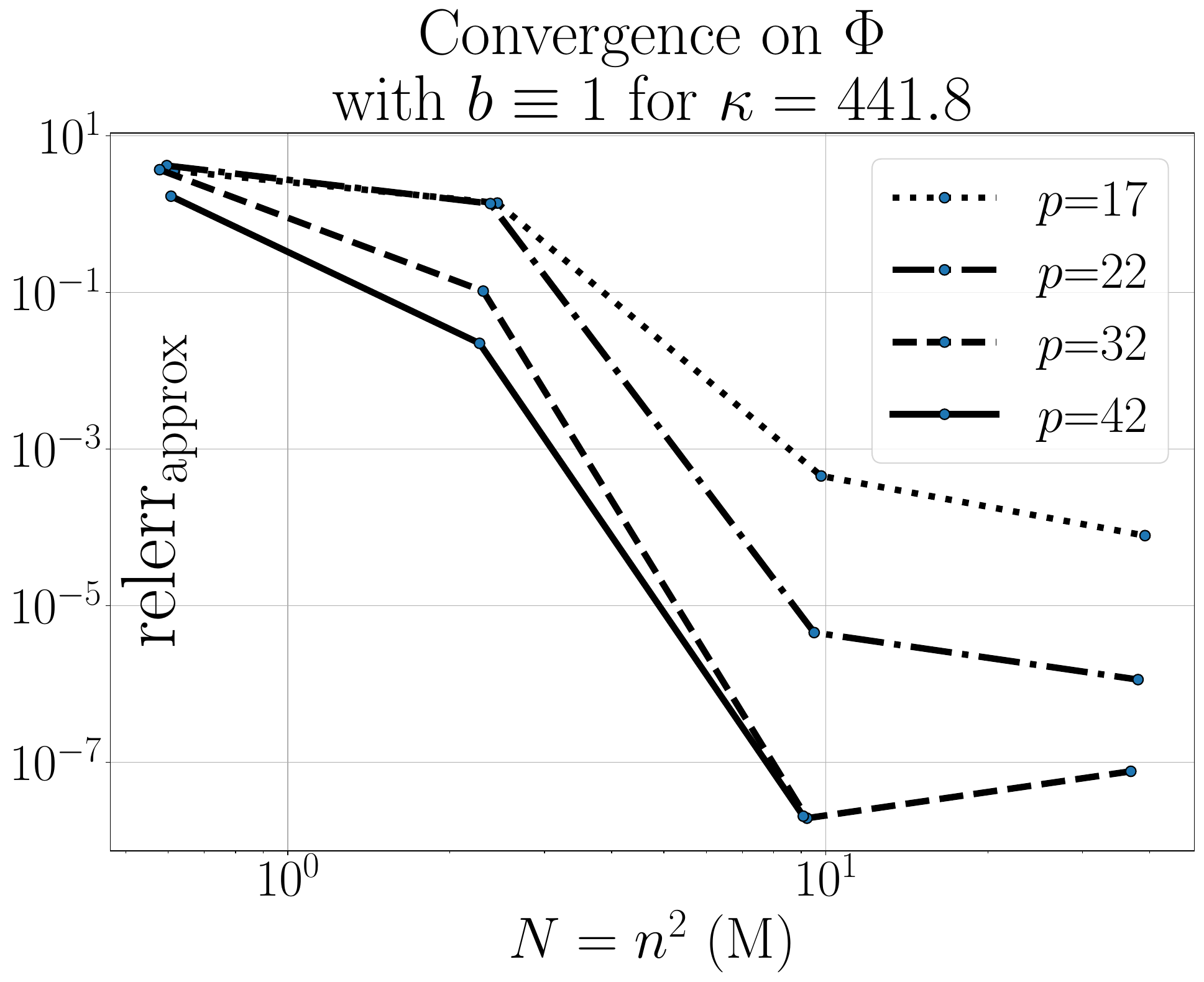}
\end{minipage}%
\hfill
\begin{minipage}{0.48\textwidth}
\caption{\small Convergence on curved domain $\Phi$ for reference solution $\mtx u_{\rm ref}$ on HPS discretization for $N$=36M with $p=42$. Choosing high orders of $p$ aids in the convergence.}
\end{minipage}
\end{subfigure}
\caption{Solutions of constant-coefficient Helmholtz problem on curved domain $\Phi$ with Dirichlet data $u \equiv 1$ on $\partial \Phi$ for various wavenumbers $\kappa$.}
\label{fig:sol_curvy_annulus}
\end{figure}

Finally, we demonstrate convergence on a curved domain $\Phi$ with a constant coefficient field $b \equiv 1$, 
where we have implemented a periodic boundary condition. The domain $\Phi$ is parametrized by the formula
\begin{equation}
    \Phi =\left\{ \left( \hat{r}(x_1,x_2) \cos \left( \hat{\theta}(x_1) \right) , \hat{r}(x_1,x_2) \sin \left( \hat{\theta}(x_1) \right) \right)\ \text{for}\ (x_1,x_2) \in [0,6] \times [0,1] \right\},
\label{eq:parametrization_curvy_annulus}
\end{equation}
where $\hat{r}(z_1,z_2) = 1 + \frac 1 5 \cos \left( \frac{15}{\pi} z_1 + z_2 \right)$
and $\hat{\theta}(z) = \frac \pi 3 z$.
By applying the chain rule, the Helmholtz operator (\ref{eq:var_helm}) on $\Phi$ takes a different form of a variable-coefficient elliptic PDE on the reference rectangle. The solutions on $\Phi$ and convergence plot are illustrated in Figure \ref{fig:sol_curvy_annulus}.

\vspace{0.5em}

\begin{remark}
When using large buffer sizes, the bound on the exact rank provided by Proposition \ref{prop:rank_property} is often pessimistic. 
For non-oscillatory problems, the actual numerical rank is very modest, even when the buffer width $b$ goes into the hundreds (cf.~Table \ref{tab:greens_hps}).
For oscillatory problems, a minimum of two points per wave-length is required, as measured at the thickest section of the slab, but the ranks are still far smaller than the strict upper bound of Prop \ref{prop:rank_property}.
\end{remark}

\section{Conclusion}
\label{sec:conc}
This paper introduces SlabLU, a sparse direct solver framework designed for solving elliptic PDEs. 
The approach decomposes the domain into a sequence of thin slabs. The degrees of freedom internal to each slab are eliminated in parallel, 
yielding a reduced matrix $\mtx T$ defined on the slab interfaces. The reduced matrix $\mtx T$
has dense sub-blocks but is rich in $\mathcal H$ and $\mathcal H^2$ matrix structure
which can be used in constructing a linear complexity direct solver for many elliptic PDEs.

One key innovation of the method is the use of randomized compression with a sparse direct solver to efficiently form $\mtx T$. 
The dense sub-blocks of $\mtx T$ have exact rank deficiencies in the off-diagonal blocks present in both the non-oscillatory and oscillatory regimes. 
The use of randomized black-box algorithms provides a purely algebraic means of efficiently forming $\mtx T$ for a variety of PDE discretizations. 

The numerical experiments presented in this paper demonstrate that SlabLU is highly effective when used in conjunction with high-order multi-domain spectral collocation schemes. 
The combination of SlabLU with high order discretization enables the rapid and accurate simulation of large-scale and challenging scattering phenomena on both rectangular and curved domains to high accuracy. The technique presented is algebraic and can readily be adapted to other standard discretization schemes such as finite element and finite volume methods
that can be partitioned into slabs.

We are currently working to further accelerate SlabLU in two regards: 
\begin{enumerate}
\item By maintaining a rank structured format for the blocks in the reduced coefficient matrix $\mtx{T}$, the memory footprint of the scheme is reduced. Moreover, such a shift gives the scheme linear complexity in the regime where the PDE is kept fixed as $N$ increases. (As opposed to the ``fixed number of points per wavelength'' scaling that is the gold standard for oscillatory problems.)
\item By replacing the sequential solve in the factorization of the reduced coefficient matrix $\mtx{T}$ by an odd-even ordering where every other block is eliminated in a hierarchical fashion, much higher parallelism can be attained.
\end{enumerate}
Both accelerations would in principle be helpful for 2D problems, but a key point of the current manuscript is that neither turned out to be necessary -- very high efficiency and essentially linear scaling is maintained up to $N \approx 10^8$ using GPU acceleration and dense 
linear algebra to factorize $\mtx T$.
In \textit{three dimensions}, however, the situation is different. Here the blocks in the reduced coefficient matrix $\mtx{T}$ hold $\sim N^{2/3}$ nodes, versus $\sim N^{1/2}$ in 2D. 
This forces us to maintain the rank structured representations of these blocks throughout the computation, in part for purposes of computational speed, but primarily to keep storage requirements from becoming excessive. 
Recent work from the authors on the randomized black-box compression and factorization 
of $\mathcal H^2$-matrices with strong admissility features useful algorithms and preliminary
experiments for compressing the sub-blocks of $\mtx T$ for thin 3D slab domains \cite{yesypenko2023randomized}. 
The extension of SlabLU to 3D is in progress, and the work will be reported at a later time.

\lsp

\noindent \textbf{Acknowledgements.} Anna would like to thank her dad, Andriy, for gifting her the RTX-3090 GPU. This version of the article has been accepted for publication, after peer review and is subject to Springer Nature’s AM terms of use, but is not the Version of Record and does not reflect post-acceptance improvements, or any corrections. The Version of Record is available online at:
\url{https://doi.org/10.1007/s10444-024-10176-x}.

\lsp

\noindent \textbf{Funding.}  The work reported was supported by the Office of Naval Research (N00014-18-1-2354), by the National Science Foundation (DMS-2313434 and DMS-1952735), and by the Department of Energy ASCR (DE-SC0022251).

% \section*{Declarations}
% \textbf{Conflict of interest} The authors have not disclosed any competing interests.

\appendix

\section{Rank Property of Thin Slabs}
\label{sec:appendix}
In this appendix, we prove Proposition \ref{prop:rank_property}, which makes
a claim on the rank structure of $\mtx T_{11}$, defined in (\ref{eq:T11_detailed}).
\rankprop*

\begin{figure}[!htb]
\centering
\begin{minipage}{0.4\textwidth}
\centering
\detailedproofpicture{28}{5}{0.20}
\end{minipage}%
\begin{minipage}{0.4\textwidth}
\captionof{figure}{To assist in the proof of Proposition \ref{prop:rank_property}, we define a partitioning of the slab interface $I_1 = J_B \cup J_F$, where $J_F = I_1 \setminus J_B$. We also partition the slab interior nodes into $I_2 := J_{\alpha} \cup J_{\beta} \cup J_{\gamma}$.}
\label{fig:partI2}
\end{minipage}
\end{figure}

Recall that $\mtx T_{11} = \mtx A_{11} - \mtx A_{12} \mtx A_{22}^{-1} \mtx A_{21}$. The proof relies on the sparsity structure of the matrices in the Schur complement. As stated in the proposition, the slab interface $I_1$ is partitioned into
indices $J_B$ and $J_F$. The proof relies on partitioning $I_2$ as well, into the indices $J_{\alpha},J_{\beta}, J_{\gamma}$ shown in Figure \ref{fig:partI2}, where $|J_{\gamma}| = 2b$.

The matrix $\mtx A_{22}$ is sparse and can be factorized as
\begin{equation}
\mtx A_{22} = \mtx L_{22} \mtx U_{22} := 
\begin{bmatrix} \mtx L_{\alpha \alpha}\\
& \mtx L_{\beta \beta}\\
\mtx L_{\gamma \alpha} & \mtx L_{\gamma \beta} & \mtx L_{\gamma \gamma}\end{bmatrix} 
\begin{bmatrix} \mtx U_{\alpha \alpha} & & \mtx U_{\alpha \gamma}\\
& \mtx U_{\beta \beta} & \mtx U_{\beta \gamma}\\
 &  & \mtx U_{\gamma \gamma}\end{bmatrix} 
 \label{eq:sparse_A22}
\end{equation}
The formula for ${\left( \mtx T_{11} \right)}_{FB}$ can be re-written as
\begin{equation}
{\left( \mtx T_{11} \right)}_{FB} = \mtx A_{FB}
 - {\left( \mtx A_{F2} \mtx U_{22}^{-1}\right)}
 {\left( \mtx L_{22}^{-1} \mtx A_{2B}\right)}
:= \mtx A_{FB} - {\mtx X_{F2}}\ {\mtx Y_{2B}}
\end{equation}
The factors $\mtx X_{F2}$ and $\mtx Y_{2B}$ have sparse
structure, due the sparsity in the factorization (\ref{eq:sparse_A22})
and the sparsity of $\mtx A_{F2}$ and $\mtx A_{2B}$.
\begin{equation}
\mtx X_{F2} = \begin{bmatrix} \mtx A_{F\alpha}& \mtx 0 & \mtx A_{F\gamma}\end{bmatrix} \mtx U_{22}^{-1},\qquad \mtx Y_{2B} 
= \mtx L_{22}^{-1} \begin{bmatrix} \mtx 0\\ \mtx A_{\beta B}\\ \mtx A_{\gamma B}\end{bmatrix}
\end{equation}
The factors $\mtx X_{F2}$ and $\mtx Y_{2B}$ have the same sparsity pattern as $\mtx A_{F2}$ and $\mtx A_{2B}$, respectively.
As a result,
\begin{equation}
{\left( \mtx T_{11} \right)}_{FB} = \mtx A_{FB} - \begin{bmatrix} \mtx X_{F\alpha}& \mtx 0 & \mtx X_{F\gamma}\end{bmatrix} \begin{bmatrix} \mtx 0\\ \mtx Y_{\beta B}\\ \mtx Y_{\gamma B}\end{bmatrix} = \underset{\text{sparse, $\mathcal O(1)$ entries}}{\mtx A_{F,B}} - 
\underset{\text{exact rank}\ 2b}{\mtx X_{F\gamma} \mtx Y_{\gamma B}}.
\end{equation}
Similar reasoning can be used to show the result for $\left( \mtx T_{11} \right)_{BF}$.

\bibliographystyle{plain}
\bibliography{main_bib}
\end{document}